\documentclass[12pt]{article}

\usepackage{lipsum}
\usepackage{amsfonts}
\usepackage{bbold}
\usepackage{graphicx}
\usepackage{epstopdf}
\usepackage{amsmath}    
\usepackage{amssymb}
\usepackage{amsthm}
\usepackage{graphicx}   
\usepackage{verbatim}   
\usepackage{color}      
\usepackage{hyperref}   
\usepackage{enumerate}
\usepackage{mathrsfs}
\usepackage{empheq}

\newtheorem{definition}{Definition}
\newtheorem{corollary}{Corollary}
\newtheorem{theorem}{Theorem}

\usepackage{algorithm}
\usepackage[noend]{algpseudocode}

\usepackage[final]{showlabels}

\newcommand{\del}{\partial}
\renewcommand{\theta}{\vartheta}
\renewcommand{\phi}{\varphi}
\newcommand{\vecc}[2]{\left ( \! \begin{array}{c}#1\\#2\\ \end{array} \! \right )}

\newcommand{\dd}{\mathrm{d}}

\newcommand{\grad}{\mathrm{grad\,}}
\renewcommand{\div}{\mathrm{div\,}}
\renewcommand{\vec}{\mathbf}

\newcommand{\vecxi}{\boldsymbol\xi}

\newcommand{\const}{\mathrm{const}}

\newcommand{\ii}{\mathbb{i}}
\newcommand{\id}{\mathbb{1}}

\newcommand{\eh}{^{\scalebox{0.5}{EH}}}
\newcommand{\ev}{^{\scalebox{0.5}{EV}}}
\newcommand{\no}{^{\scalebox{0.5}{N}}}
\newcommand{\avg}{}

\newcommand{\doubleImage}[6]{%
\begin{figure}[!bt]%
\begin{minipage}{0.45\textwidth}
\centering
\includegraphics[width=\textwidth]{#1}
\caption{#2}%
\label{#3}%
\end{minipage}
\hfill
\begin{minipage}{0.45\textwidth}
\centering
\includegraphics[width=\textwidth]{#4}
\caption{#5}%
\label{#6}%
\end{minipage}
\end{figure}
}

\numberwithin{equation}{section}

\begin{document}

\begin{center} \Large
The active flux scheme on Cartesian grids and its low Mach number limit\footnote{We thank Yifan Bai of the University of Michigan for providing the results shown in Figure \ref{fig:unstructuredgresho}. WB acknowledges support of the German National Academic Foundation and the support of the German Academic Exchange Service (DAAD) with funds from the German Federal Ministry of Education and Research (BMBF) and the European Union (FP7-PEOPLE-2013-COFUND -- grant agreement no. 605728).}

\vspace{1cm}

\date{}
\normalsize

Wasilij Barsukow\footnote{Institute of Mathematics, Zurich University, 8057 Zurich, Switzerland \emph{and} Institute of Mathematics, Wuerzburg University, Emil-Fischer-Strasse 40, 97074 Wuerzburg, Germany, \texttt{wasilij.barsukow@math.uzh.ch} } \hspace{2cm}
        Jonathan Hohm\footnote{Institute of Mathematics, Wuerzburg University, Emil-Fischer-Strasse 40, 97074 Wuerzburg, Germany} \\
        Christian Klingenberg\footnote{Institute of Mathematics, Wuerzburg University, Emil-Fischer-Strasse 40, 97074 Wuerzburg, Germany} \hspace{2cm}
        Philip L. Roe\footnote{Department of Aerospace Engineering, University of Michigan, 1320 Beal Avenue, Ann Arbor, MI 48109, USA.}
\end{center}




\begin{abstract}

Finite volume schemes for hyperbolic conservation laws require a numerical intercell flux. In one spatial dimension the numerical flux can be successfully obtained by solving (exactly or approximately) Riemann problems that are introduced at cell interfaces. This is more challenging in multiple spatial dimensions. The \emph{active flux scheme} is a finite volume scheme that considers continuous reconstructions instead. The intercell flux is obtained using additional degrees of freedom distributed along the {cell boundary}. For their time evolution an exact evolution operator is employed, which naturally ensures the correct direction of information propagation and provides stability. This paper presents an implementation of active flux for the acoustic equations on two-dimensional Cartesian grids and demonstrates its ability to simulate discontinuous solutions with an explicit time stepping in a stable manner. Additionally, it is shown that the active flux scheme for linear acoustics is low Mach number compliant without the need for any fix.


\end{abstract}

\section{Introduction}

Conservation laws, such as the Euler equations of ideal hydrodynamics, express conservation of density, momentum and energy of a compressible fluid. Numerical methods need to reflect the conservation property, as by the Lax-Wendroff theorem only conservative numerical schemes are able to converge to a (weak) solution of the equations. Cell-based methods consider a partition of space into computational cells. A natural class of numerical schemes for conservation laws are finite volume schemes. They interpret the discrete degree of freedom in a computational cell as the average of the dependent variable. Numerical fluxes assigned to the {cell boundaries} automatically ensure conservativity of the method.

In \emph{one spatial dimension} the numerical flux can be successfully obtained from a so-called Riemann solver. The dependent variables are reconstructed in a piecewise constant manner allowing for jumps at the locations of cell interfaces. The Riemann problems that arise at cell interfaces are then solved (see e.g. \cite{godunov59difference}). It is often possible to solve the Riemann problem exactly, but an approximate solution can lead to a similarly accurate scheme while requiring less computations. Therefore approximate Riemann solvers are very popular (e.g. \cite{bouchut04,toro09}). Among others, relaxation (e.g. \cite{jin95}) can be a successful strategy to construct them.

In \emph{multiple spatial dimensions} one faces additional challenges. First of all, even approximate Riemann solvers are complicated (see e.g. {\cite{zheng12,barsukow17}} for exact solutions to multi-dimensional Riemann problems and e.g. \cite{colella90,balsara12} for approximate solvers). This explains why dimensionally split methods are so popular: they consider the Riemann problems arising at the edges of the computational cell independently. The multi-dimensional problem is thus replaced by several one-dimensional ones. 

Additionally, the solutions to Euler equations in multiple spatial dimensions exhibit a number of phenomena absent in one dimension. Particularly prominent features are vortices, that strongly dominate realistic multi-dimen\-sio\-nal flows. In the limit of low Mach numbers, the solutions of compressible Euler equations tend to those of incompressible Euler equations. This is a multi-dimensional feature as well, because incompressible flows are trivial in one spatial dimension. Many finite volume methods suffer from artefacts when applied e.g. to the regime of low Mach number (see e.g. \cite{guillard99}). While Riemann solvers seem to yield very good results in one spatial dimension it is unclear whether this is true in multiple spatial dimensions, or whether numerical methods based on different concepts are more adequate. It has been, for instance, noticed that the failure in the regime of low Mach number is observed even when using \emph{exact} multi-dimensional Riemann solvers (\cite{guillard04,barsukow17}).

It thus seems relevant to investigate alternative ways how multi-dimensional numerical methods can be constructed. Recently, such a method has been proposed by \cite{eymann13}, as an extension of a method from \cite{vanleer79} for linear advection. It has been given the name \emph{active flux}. It involves a continuous reconstruction, together with pointwise degrees of freedom along the cell boundary. It thus does not require the solution of any Riemann problem. Continuous reconstructions are found not to stand in the way of computing discontinuous solutions. Indeed, even schemes that employ reconstructions with jumps across cell interfaces are mostly unable to resolve a shock wave sharply. They resolve a shock wave by a succession of jumps. Instead, a method with a continuous reconstruction resolves the {discontinuity} by a continuous function with a sharp gradient (see Figure \ref{fig:afrecondetail}). The details of the method are discussed in section \ref{structureAF} below.

\begin{figure}[h]
 \centering
 \includegraphics[width=0.45\textwidth]{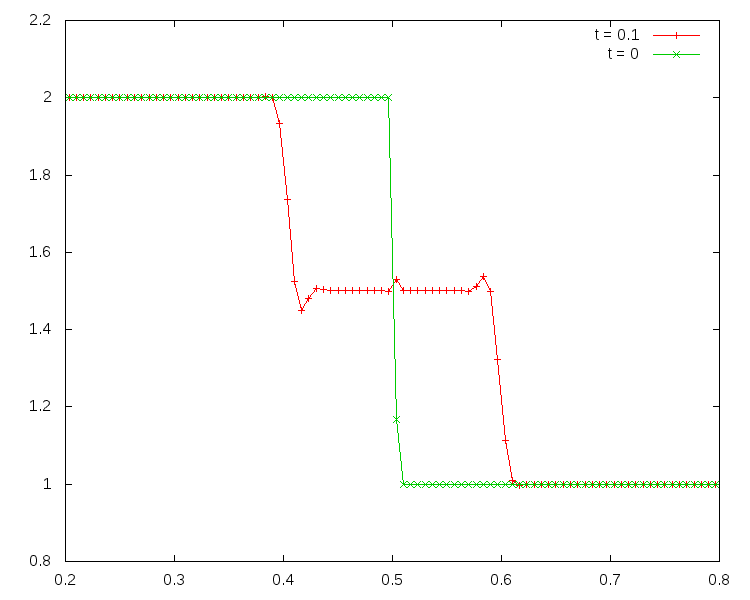}\hfill\includegraphics[width=0.45\textwidth]{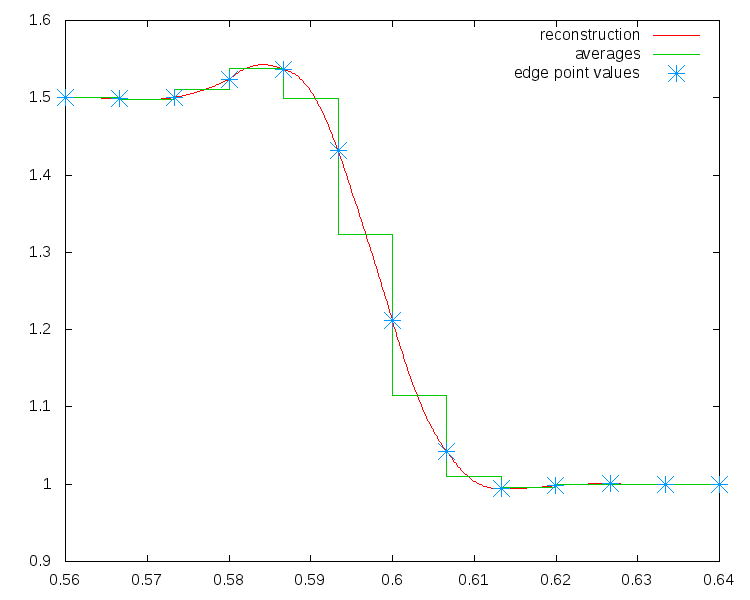}
 \caption{Illustration of the continuous reconstruction in the context of discontinuous solutions. The figures show the solution to a Riemann problem of linear acoustics with the active flux scheme in one spatial dimension on a grid of 100 cells. The exact solution consists of two {discontinuities} moving with speeds $\pm1$. \textit{Left}: Initial data and the numerical solution at time $t=0.1$. The overshoots are due to the scheme being higher order, but no limiting employed. \textit{Right}: Close-up of the solution at the location of the right {discontinuity}. The solid line is the continuous reconstruction; additionally the averages and the point values at cell boundaries are shown.}
 \label{fig:afrecondetail}
\end{figure}

This paper aims at contributing to a deeper understanding of the active flux method. First, a very general formulation of the method is given, that is independent of the grid or the equations (section \ref{structureAF}). It thus highlights the most important features that distinguish the method from conventional finite volume schemes. Active flux in multiple spatial dimensions has so far been introduced for triangular computational cells. For certain applications triangles might be too complicated a mesh. This paper for the first time presents an implementation on {two-dimensional} Cartesian grids (section \ref{sec:cartesianaf}). 

One of the objectives behind the development of the active flux method is to improve multi-dimensional simulations by including more multi-dimensional information in the scheme. {Following \cite{eymann13}, in this paper the active flux method is applied to the equations of linear acoustics.} Linear acoustics is obtained when the Euler equations are linearized around a constant background with no velocity. In multiple spatial dimensions, linear acoustics cannot be reduced to some kind of multi-dimensional advection. This makes it a valuable system in order to investigate the performance of numerical methods in multiple spatial dimensions (as e.g. also in \cite{dellacherierieper10}). 

The active flux scheme pursues a strategy alternative to discontinuous reconstructions. The fact that the method does not use a Riemann solver raises the question whether maybe it is better suited for the low Mach number limit than Riemann solvers. Linear acoustics exhibits a low Mach number limit similar to that of the Euler equations and there exists an efficient methodology to study the behaviour of numerical schemes in this limit on Cartesian grids (introduced in \cite{barsukow17a,barsukow17lilleproceeding}). Here, this methodology is extended to take into account the {additional degrees of freedom along the cell boundary} and is applied to the active flux scheme in section \ref{statPresAF}. It is shown that indeed, active flux is able to resolve the low Mach number limit of linear acoustics. 

{The presence of point values along the boundary makes it easy to evaluate the fluxes necessary for the update of the cell average. These point values, however, also require a time evolution. Here, following \cite{eymann13}, the exact solution operator of linear acoustics is applied onto reconstructed initial data. However, this represents only one possible choice of an evolution operator for the point values and other evolutions are possible. The presentation of the general philosophy of active flux in section \ref{structureAF} therefore does not insist on a particular way of evolving the point values. This becomes particularly important for nonlinear equations without an available exact evolution operator. Recent work on approximate evolution operators for nonlinear equations is \cite{kerkmann18}.}

\section{Acoustic equations}\label{chapter2dAcoustic}

This paper considers numerical methods for $n \times n$ systems of conservation laws in $d=2$ spatial dimensions:
\begin{align}
 \del_t q + \nabla \cdot \vec f(q) &= 0  \label{eq:conslaw}\\
 q &: \mathbb R^+_0 \times \mathbb R^d \to \mathbb R^n\\
 \vec f = (f^x, f^y) \qquad 
 f^x, f^y &: \mathbb R^n \to \mathbb R^n
\end{align}

Only vectors with $d$ components are set in boldface symbols. Additionally, indices are never denoting a derivative.

The most prominent example of a system of conservation laws are the Euler equations 
\begin{align}
\begin{split}
		\del_t \rho+ \nabla\cdot(\rho \vec v)&=0 \\
		\del_t (\rho\vec{v}) +\nabla\cdot(\rho \vec{v}\otimes\vec{v} + p \id)&=0
\end{split}
\label{eulerEqKons}
\end{align}
with velocity $\vec{v}=(u,v)^T$, pressure $p$ and density $\rho$.

As a stepping stone towards a detailed understanding of numerical methods for the Euler equations in multiple spatial dimensions the \emph{acoustic equations} are studied first. They are obtained as a linearization of the Euler equations around the state of constant density and pressure and vanishing velocity (see also e.g. \cite{barsukow17}). Their symmetrized version reads

\begin{align}
	\begin{split}
		\del_t p+ c\nabla \cdot \vec{v}&=0, \\
		\del_t \vec{v}+ c\nabla p&=0.
	\end{split}
	\label{acousticEq}
	\end{align}
with $c > 0$ the speed of sound of the background state. This system is strongly hyperbolic in two spatial dimensions with eigenvalues $\{ \pm c, 0 \}$. It has to be augmented by initial data 
\begin{align}
p(0, \vec x) &= p_0(\vec x) & \vec v(0, \vec x) &= \vec v_0(\vec x) \label{eq:initialdataacoustic}
\end{align}
and possibly boundary conditions. 

The curl of the velocity is called \emph{vorticity}. {In the context of Euler equations, there is great interest in methods that are able to compute reliably flows that contain vorticity, so that the appearance of vortical structures is caused by physics and affected as little as possible by numerical error.

For the acoustic equations vorticity is stationary}:
\begin{align}
 \del_t (\nabla \times \vec v) &= 0 \label{eq:stationaryvorticity}
\end{align}
Thus, for acoustics, an interesting class of methods are those that keep a discrete version of the vorticity exactly stationary (\emph{vorticity preserving}), see e.g. \cite{morton01,jeltsch06,mishra09preprint,barsukow17a}.

The low Mach number limit for the equations of linear acoustics is the limit $\epsilon \to 0$ of
\begin{align}
	\begin{split}
		\del_t p+ \frac{c}{\epsilon}\nabla \cdot \vec{v}&=0, \\
		\del_t \vec{v}+ \frac{c}{\epsilon}\nabla p&=0.
	\end{split}
	\label{acousticEqEpsilon}
	\end{align}
{It is obtained by analogy from the low Mach number limit of the Euler equations. There, the pressure gradient $\nabla p$ is rescaled as $\frac{\nabla p}{\epsilon^2}$. The same is done here for acoustics, which after symmetrization yields \eqref{acousticEqEpsilon}. The interpretation of the Mach number $\epsilon$ as the ratio between advective and acoustic speeds in the setting of the Euler equations, however, is not applicable to the acoustic equations. Nevertheless this limit is useful in the design and understanding of low Mach number compliant schemes.}

The low Mach number limit in \eqref{acousticEqEpsilon} is the same as the limit $t \to \infty$ {of the initial value problem \eqref{acousticEq}--\eqref{eq:initialdataacoustic}}. With appropriate boundary conditions, the long time solution is governed by the stationary states (for more details see \cite{barsukow17a}).

Note that system \eqref{acousticEq} can be reduced to a wave equation for $p$, and a \emph{vector wave equation}
\begin{align}
 \del_t^2 \vec v - c^2 \nabla (\nabla \cdot \vec v) &= 0
\end{align}
for $\vec v$. This latter, however, cannot generally be reduced to several scalar wave equations for the components of $\vec v$, because $\nabla(\nabla \cdot \vec v) - \Delta \vec v = \nabla \times ( \nabla \times \vec v)$.

The numerical method will make use of the exact solution to the initial value problem \eqref{acousticEq}--\eqref{eq:initialdataacoustic} for non-differentiable initial data. The required (distributional) solution operator for initial data of such low regularity has been derived in detail in \cite{barsukow17}. The solution formulae contain spherical means, that also appear in the solution of the scalar wave equation (\cite{john78}). It is helpful to first consider three spatial dimensions:

\begin{definition}[Spherical mean]
The \emph{spherical mean} $M\left[f\right](\vec{x}, r)$ of an integrable function $f$ that depends on $\vec x = (x, y, z) \in \mathbb R^3$ is given by
\begin{align}
\begin{split}
	M \left[f\right](\vec{x}, r)&=\frac{1}{4\pi} \oint_{S^2} \dd \vec y \, f(\vec x + r \vec y)
	=\frac{1}{4 \pi}\int\displaylimits_0^{2\pi} \dd\phi \int\displaylimits_0^\pi \dd\theta \, \sin \theta f\left(\vec{x}+r\cdot \vec{n}\right) 
\end{split}
\label{SM3D}
\end{align}
with the outward normal vector given by
\[
\boldsymbol{n} = 
\begin{pmatrix}
\sin\theta\cos\phi \\
\sin\theta\sin\phi \\
\cos\theta
\end{pmatrix}
\]
\end{definition}

{This paper concentrates on the two-dimensional case, though. If $f$ does not depend on $z$, then, with $s := r \sin \theta$, it is possible to rewrite equation \ref{SM3D} as
\begin{align}
	M^{\scalebox{0.5}{$2D$}}\left[f\right](x, y, r)=\frac{1}{2 \pi r}\int\displaylimits_0^{2\pi} \dd \phi \int\displaylimits_0^r \dd s f\left(x+s \cos\phi,y+s \sin\phi\right)\frac{s}{\sqrt{r^2-s^2}} 
\label{SM2D}
\end{align}

In both cases, if $f$ is polynomial, the spherical means can be evaluated analytically. Indeed, by shifting the point of integration one is left with 
\begin{align}
	M^{\scalebox{0.5}{$2D$}}\left[x^py^q\right](0,0, r) &= \frac{1}{2 \pi r} \int\displaylimits_0^{2\pi} \dd \phi  \cos^p\phi \sin^q\phi  \int\displaylimits_0^r \dd s \frac{s^{1+p+q}}{\sqrt{r^2-s^2}}
\end{align}
Both integrals can easily be evaluated analytically because the angular integration bounds here are multiples of $\frac{\pi}{2}$. The second integral is
\begin{align}
  \int\displaylimits_0^r \dd s \frac{s^m}{\sqrt{r^2-s^2}} 
  &= \begin{cases} \displaystyle r^{2m'} \frac{\pi }{ 4^{m'} } { 2 m'-1 \choose m'}  & m = 2m', m' \in \mathbb N \\ \\
  \displaystyle r^{2m'+1} \frac{ (m')!^2 4^{m'}}{(2m'+1)!  } &m = 2 m' + 1, m' \in \mathbb N   \end{cases} \label{eq:sphmeanpoly}
\end{align}

}

In \cite{eymann13} for the initial value problem \eqref{acousticEq}--\eqref{eq:initialdataacoustic} the following solution formula appears:
\begin{align}
		p(t,\vec{x})&=p_{0}(\vec{x})+\int\displaylimits_0^{ct} \dd r \, r\cdot M\big[\div \grad p_{0} \big](\vec x, r)-ct\cdot M\big[\div \vec{v}_{0} \big](\vec x, ct) \label{solGradDivP}\\
		\vec{v}(t,\vec{x})&=\vec{v}_{0}(\vec{x})+\int\displaylimits_0^{ct}\dd r \, r\cdot M\big [\grad \div \vec{v}_{0} \big ](\vec x,r) -ct\cdot M\big[\grad p_{0}\big](\vec x, ct) \label{solGradDivU}
\end{align}
{These formulas are derived from the classical Poisson formula for the scalar wave equation in second-order form \cite{whitham11}, which was used in \cite{alpert00,hagstrom15} to derive high-order finite-difference methods for the scalar problem. For the acoustic equations in first-order form, and for the linearized Euler equations, \eqref{solGradDivP} was applied to the pressure, and a similar equation to the velocity components, and employed in \cite{eymann13} to create an active flux scheme, with vorticity treated as a source term. The full extension to vortical flows \eqref{solGradDivP}--\eqref{solGradDivU} appears in \cite{fan15}, together with investigations of the active flux method on unstructured grids, including boundary conditions and nonlinearity. Finite-volume schemes based on the full extension are reported in \cite{franck17}}.

For the numerical method it will be necessary to consider initial data that are continuous, but possess discontinuous first derivatives, such that terms like $\nabla \cdot \nabla p_0$ require clarification. The interpretation of this formula in the sense of distributions has been achieved in \cite{barsukow17}. At the same time it has been shown that the formula can be rewritten as follows:

\begin{theorem}[Solution operator]
The solution to the initial value problem \eqref{acousticEq}--\eqref{eq:initialdataacoustic} is given by
\begin{align}
	p(t,\vec{x})=&\,\del_r\Big(r\cdot M\left[p_{0}\right](\vec x, r)\Big)\Big|_{r=ct}-\frac{1}{ct}\del_r\Big(r^2 M\left[\vec{v}_{0}\cdot\vec{n}\right](\vec x, r)\Big)\Big|_{r=ct} \label{eq:solp}\\
	\begin{split}
\vec{v}(t,\vec{x})=&\,\vec{v}_{0}(\vec x)-\frac{1}{ct}\del_r\Big(r^2 M\left[p_{0}\vec{n}\right](\vec x, r)\Big)\Big|_{r=ct}\\
&+\int\limits_0^{ct} \frac{1}{r}\del_r\left(\frac{1}{r}\del_r\left(r^3 M\left[\left(\vec{v}_{0}\cdot\vec{n}\right)\vec{n}\right](\vec x, r)\right)-r M\left[\vec{v}_{0}\right](\vec x, r)\right)\dd r \end{split}\label{eq:solu}
\end{align}
The derivatives are to be interpreted in the sense of distributions, if necessary.
\end{theorem}
For the proof of the theorem, and more details on its distributional version, see \cite{barsukow17}. Below it is shown that for the reconstructions considered here all the derivatives in \eqref{eq:solp}--\eqref{eq:solu} exist in the strong sense and it is not necessary to consider anything distributionally. When trying to use formula \eqref{solGradDivP}--\eqref{solGradDivU}, one would be forced to interpret them in the sense of distributions in order to compute the second derivatives correctly. Nothing seems to be gained when using \eqref{eq:solp}--\eqref{eq:solu}, because they also contain second derivatives. However this time, the derivatives are all expressed with respect to $r$ and the initial data in the particular setup of this paper will turn out to have continuous derivatives with respect to $r$ (section \ref{chapterRectangularUpdate}). Thus when using formulae \eqref{eq:solp}--\eqref{eq:solu} it is possible to avoid dealing with distributional solutions, which is a substantial advantage.

It should be noted, that the exact solution may be expressed in various representations, which differ by the way they can be used. In this respect the above formulae are very different from those obtained, for example, in \cite{ostkamp97,lukacova04a} using bicharacteristics. There, analytical relations are derived, which connect the solution at time $t > 0$ with the data at initial time via a so called \emph{mantle integral}. This latter involves the solution at all intermediate times. Therefore, the knowledge of the data at initial time does not allow to immediately compute the solution at a later time \emph{exactly}. With the above formulae \eqref{solGradDivP}--\eqref{solGradDivU} or \eqref{eq:solp}--\eqref{eq:solu}, on other hand, the solution at time $t$ is given as a functional of the data at initial time only and thus can be determined immediately.

\section{General structure of the active flux method}\label{structureAF}

This section outlines the idea of the active flux method on general polygonal (unstructured) meshes and without specifying the hyperbolic system of PDEs that is to be solved. Several references to implementations of active flux for particular equations that appear in the literature are given in section \ref{sec:standardaf}. In section \ref{sec:cartesianaf} finally, Cartesian meshes with rectangular cells are introduced, on which the active flux method is used to solve the acoustic equations \eqref{acousticEq}. 

The active flux method is an extension of the finite volume method. Therefore finite volume methods are first reviewed in section \ref{sec:fv}, before the active flux method is introduced in section \ref{sec:af}.

\subsection{Finite volume methods} \label{sec:fv}

For the finite volume method the discrete degree of freedom $q_\mathcal C$ is interpreted as the average of $q$ over the polygonal numerical cell $\mathcal C \subset \mathbb R^d$. For Cartesian grids, also the notation $q_{\mathcal C_{ij}} =: q_{ij}$ is standard. However, our presentation of the method for the moment does not depend on the nature of the grid. In order to construct the method one first integrates the conservation law \eqref{eq:conslaw} over the cell and applies Gauss' law. This gives rise to fluxes through the cell boundary:
\begin{align}
 \del_t \int_{\mathcal C} \dd \vec x \, q + \int_{\del \mathcal C} \dd \vec x \, \vec n \cdot \vec f(q) &= 0 \label{eq:fvexact}
\end{align}
Here $\vec n$ is the outward normal on the {cell boundary} $\del \mathcal C$ of cell $\mathcal C$.

The finite volume method now replaces the exact expression \eqref{eq:fvexact} by a (yet undetermined) numerical approximation $\bar f_e$ of the flux through an edge $e \subset \mathcal \del C$:
\begin{align}
 \del_t q_\mathcal C + \frac{1}{|\mathcal C|} \sum_{e \subset \del\mathcal C} |e| \bar f_e &= 0 \label{eq:fv}
\end{align}
Here $|e|$ is the length of the edge $e$ and $|\mathcal C|$ the area of the cell $\mathcal C$.

One way to obtain a numerical flux is by using a Riemann solver: One considers the initial value problem \eqref{eq:conslaw} with initial data $q_{\text{recon}}: \mathbb R^d \to \mathbb R^n$ that are piecewise constant
\begin{align}
 q_{\text{recon}}(\vec x) = q_\mathcal C \quad \text{if} \quad \vec x \in \mathcal C 
\end{align}
This step is called \emph{reconstruction}. The resulting Riemann problems across the cell boundaries are solved exactly or approximately in time. The numerical flux is then taken to be the time average of the flux of this solution through the {cell boundary}. 

Alternatively, one might consider an initial value problem given by more complicated reconstructions inside the cells, which is referred to as \emph{generalized Riemann problem} {(see e.g. \cite{glimm84,bourgeade89,benartzi03,titarev02,toro02}} and many others). This leads to schemes that are higher order discretizations of the original PDE. In this case the reconstruction $q_{\text{recon}}$ in cell $\mathcal C$ does not only depend on $q_\mathcal C$, but also on the values in neighbouring cells. Conservation still requires the average $\frac{1}{|\mathcal C|} \int_\mathcal C \dd \vec x \, q_{\text{recon}}(\vec x)$ of the reconstruction over any cell $\mathcal C$ to match the discrete degree of freedom $q_\mathcal C$ in this cell.

\subsection{Active flux method} \label{sec:af}

Now the active flux scheme is introduced. Its presentation in this section again is given without specifying the type of grid. In section \ref{sec:standardaf} a brief overview of an implementation on one-dimensional grids and on triangles is given, before section \ref{sec:cartesianaf} presents the novel implementation on Cartesian grids.

The active flux method differs from usual finite volume methods in how the numerical flux is obtained. Additional degrees of freedom\footnote{Recall that throughout the paper indices never denote a derivative.} $q_{\vec p}$, $\vec p \in \del \mathcal C$ are introduced, which are given the interpretation of point values and are distributed at a finite number of locations $\vec p$ along the cell boundary $\del \mathcal C$. This immediately allows to use them in a quadrature formula in order to approximate the flux through the {cell boundary}. The update formula for the cell average $q_\mathcal C$ remains equation \eqref{eq:fv} just as for usual finite volume methods. 

The additional degrees of freedom require {an update procedure which shall be detailed next}. Here, a reconstruction $q_{\text{recon}}:\mathbb R^d \to \mathbb R^n$ is used as initial data. {The resulting initial value problem \eqref{eq:conslaw} is then solved (exactly or approximately)} at the location $\vec p \in \del \mathcal C$ of the {point value at the cell boundary}. {In the implementation of active flux for the acoustic equations below, an exact solution of the initial value problem is used (equations \eqref{eq:solp}--\eqref{eq:solu}).}

The reconstruction $q_{\text{recon}}$ has to fulfill both
\begin{align}
 \frac{1}{|\mathcal C|} \int_\mathcal C \dd \vec x \, q_{\text{recon}}(\vec x) &= q_\mathcal C \label{eq:conservation}
\end{align}
and
\begin{align}
 q_{\text{recon}}(\vec p) &= q_{\vec p} \qquad \text{at finitely many } \vec p \in \del \mathcal C \label{eq:boundaryinterpol}
\end{align}
The latter condition implies that the reconstruction is continuous across the cell boundary at least in all points $\vec p$. In fact, in this paper the reconstructions will be continuous across all of $\del \mathcal C$. The increased number of conditions makes the appearance of high order interpolation polynomials natural, which make the scheme high order in space. High order temporal accuracy is obtained by using one or several substeps during the time integration (as discussed in section \ref{chapterRectangularFlux}).

\begin{figure}[b!]%
\centering
\includegraphics[width=.7\columnwidth]{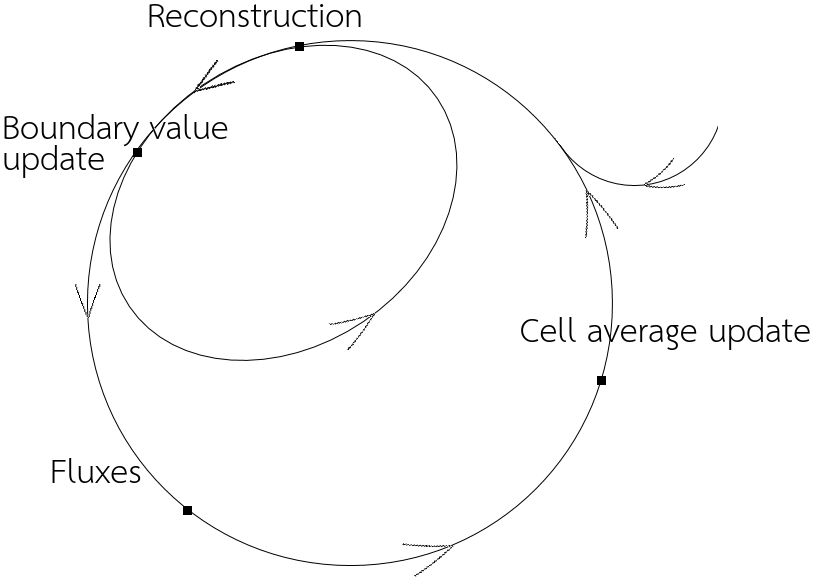}
\caption{The active flux algorithm.}%
\label{flowchart1D}%
\end{figure}

To summarize, the active flux method promotes the solution along the cell boundary (and thus in a sense the cell boundary flux) to the status of an independent degree of freedom (hence the name \emph{active flux}). In usual finite volume methods on the other hand the cell boundary flux is a derived quantity. Another difference to usual finite volume methods is that the reconstruction is used to evolve the {point values along the cell boundary} (which are then used to compute the flux), rather than to obtain a flux directly by e.g. solving a Riemann problem. The algorithmic structure of the active flux method is shown in Figure \ref{flowchart1D} and in algorithm \ref{AFalgorithm}.

\begin{algorithm}
\caption{Active flux algorithm}\label{AFalgorithm}
\begin{algorithmic}[1]
\Procedure{ActiveFlux}{$numberOfTimeSteps$}
\While{$t<t_{final}$}
\For{$cells$}
\State build reconstruction for the cell
\EndFor\label{reconstuctFor}
\For{$cell boundaries$}
\State time evolution of the point values at {cell boundaries}
\EndFor\label{updateFor}
\For{$cell boundaries$}
\State calculate numerical flux at the {cell boundary}
\EndFor\label{numericalFluxFor}
\For{$cells$}
\State calculate new average for the cell
\EndFor\label{fvmFor}
\State $t$\verb-+=-$\Delta t$
\EndWhile
\EndProcedure
\end{algorithmic}
\end{algorithm}

\subsection{The standard active flux method} \label{sec:standardaf}

The active flux method as formulated in section \ref{sec:af} is very general. It is not restricted to a particular grid, nor to a particular form of the conservation law \eqref{eq:conslaw}. At the same time it leaves a lot of freedom concerning its design, as neither the precise locations of the {point values along the cell boundary} nor the reconstruction have been specified in section \ref{sec:af}.

Two particular designs appear in the literature and shall be briefly discussed here for the sake of an overview. The first goes back to \cite{vanleer79} and considers a one-dimensional grid. The cell boundaries are points, and an additional degree of freedom is located at each one. The reconstruction is piecewise parabolic, which is the lowest degree polynomial that fulfills conditions \eqref{eq:conservation} and \eqref{eq:boundaryinterpol} in this setting. The reconstruction is continuous across the cell boundaries due to \eqref{eq:boundaryinterpol}, but has discontinuous derivative there. \cite{vanleer79} discusses the application to linear advection; Burgers' equation has been considered in \cite{eymann11a,roe17a}, nonlinear hyperbolic systems in \cite{eymann11}.

An extension to two-dimensional triangular grids has been introduced in \cite{eymann13,eymann13a}. Inspired by P2 interpolations the {point values} are located at vertices and edge midpoints of the triangles. \eqref{eq:conservation} and \eqref{eq:boundaryinterpol} are 7 equations in this case. The reconstruction is taken to be a subset of the biparabolic polynomials (see \cite{eymann13}, Table 2 for details). In \cite{eymann13} the method is used to numerically solve both linear advection and linear acoustics. Extensions to nonlinear problems are discussed in \cite{fan17,maeng17}.

In both cases the reconstructions are piecewise parabolic, thus yielding third order accuracy in space. In order to reach the same accuracy for the temporal discretization, the {point values along the cell boundary} are evolved with half time steps {from the same data}. The quadrature formula for the numerical flux then is a space-time Simpson rule.

\section{Active flux on a Cartesian mesh for the acoustic equations}  \label{sec:cartesianaf}
In this chapter, an implementation of the active flux scheme on rectangular meshes is considered for the first time for the acoustic equations. Acoustic equations in multiple spatial dimensions are an important stepping stone in the development of the active flux method for the Euler equations. At the same time Cartesian grids are easy to implement and are used in a variety of applications. 

A Cartesian mesh consists of rectangular cells $\mathcal C_{ij}$ with width $\Delta x$ and height $\Delta y$ 
\begin{align}
\mathcal C_{ij} = \left[ \left(i-\frac12\right)\Delta x,\left(i+\frac12\right)\Delta x \right ] \times \left [ \left(j-\frac12\right)\Delta y, \left(j+\frac12\right)\Delta y\right ] \subset \mathbb R^2
\end{align}
indexed by $(i,j) \in \mathbb Z^2$. Therefore $|\mathcal C| = \Delta x \Delta y$. The cell average $q_{\mathcal C_{ij}}$ is denoted by $q_{ij}$.

The algorithm consists of the four steps mentioned in section \ref{sec:af}. The distribution of the {point values along the cell boundary} is presented in section \ref{chapterboundarydeg}. The reconstruction on rectangular cells is discussed in section \ref{chapterRectangularReconstruction}. It is used as initial data in order to advance the {point values at the cell boundaries} forward in time. The corresponding solution operator for the equations of linear acoustics is presented in section \ref{chapterRectangularUpdate}. A quadrature rule is then applied to the {point values along the cell boundaries} in order to determine the numerical flux (section \ref{chapterRectangularFlux}), and the fluxes are used to update the cell average.

\subsection{{Distribution of point values along the cell boundary}}\label{chapterboundarydeg}

In this and the next section, the reconstruction is adapted to a rectangular cell. As on triangular grids mentioned in section \ref{sec:standardaf}, the {point values} are taken to be located at the corners of the cell (\emph{node values}) and at edge centers (\emph{edge values}). These locations are indexed as in Figure \ref{EdgeNodeRectangle}. Each cell $\mathcal C_{ij}$ has eight {point values along the cell boundary}: four node and four edge values. These degrees of freedom are correspondingly denoted by $q_{m,ij}$, $m=1,...,8$. They are shared by the adjacent cells, i.e. for example $q_{4, ij} = q_{8, i+1,j}$.

\doubleImage{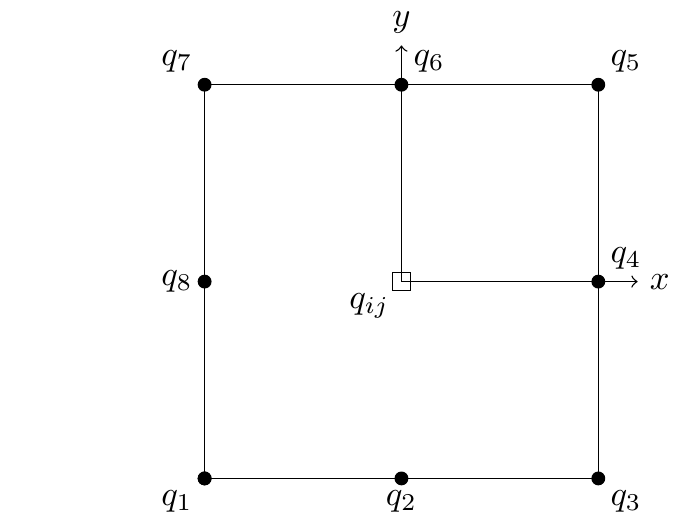}{The {point values} $q_{1,ij}, ..., q_{8,ij}$ {located along the cell boundary} and the cell average $q_{ij}$ for a rectangular cell $\mathcal C_{ij}$ with width $\Delta x$ and height $\Delta y$.}{EdgeNodeRectangle}{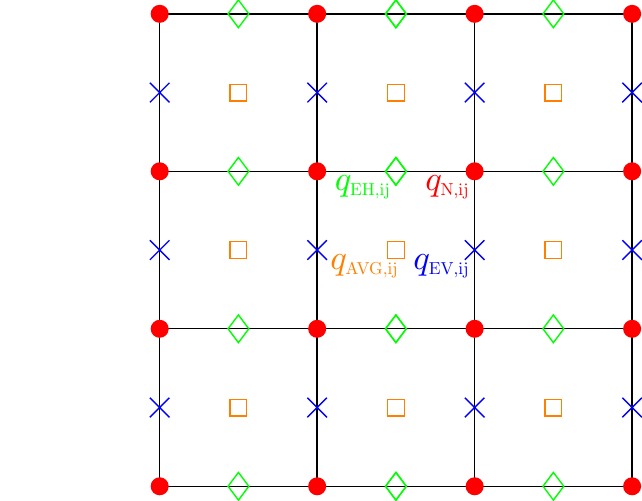}{The grid for the nodes (red circles), the vertical edges (blue crosses), the horizontal edges (green diamonds) and the averages (orange squares).}{stationarityGrids}

For usual finite volume methods, the cell average is the only degree of freedom per cell. Now additionally to the cell average there are 8 degrees of freedom distributed {as point values} along the cell boundary: There are four node values, each shared by four cells. Similarly, there are four edge values, which are shared by two cells each. Therefore \emph{per cell} one is left with
\begin{itemize}
 \item 1 cell average
 \item 1 node value
 \item 2 edge values (horizontal and vertical)
\end{itemize}
In Figure \ref{stationarityGrids} these degrees of freedom are shown in different colors. Indeed, each of them forms a lattice with spacings $\Delta x$ and $\Delta y$ in $x$ and $y$ direction, respectively. When solving linear acoustics with the active flux method, there are thus 12 variables per cell, i.e. 4 degrees of freedom per cell with 3 variables per degree of freedom.

\subsection{Reconstruction}\label{chapterRectangularReconstruction}

The reconstruction is subject to the constraints \eqref{eq:conservation} and \eqref{eq:boundaryinterpol} from section \ref{structureAF}: it has to be exact in the eight {point values along the cell boundary} and the average of the reconstruction over the cell has to yield the cell average. These are nine conditions. It is helpful to define the reconstruction $q_{\text{recon},ij}$ in any cell to refer to the cell midpoint as $\vec x = 0$. Therefore denote by
\begin{align}
q_{\text{recon},ij} : \left[-\frac{\Delta x}{2}, \frac{\Delta x}{2}\right] \times \left[-\frac{\Delta y}{2}, \frac{\Delta y}{2}\right] \to \mathbb R^n \label{eq:reconallgrid}
\end{align}
the reconstruction of any quantity $q$ in cell $\mathcal C_{ij}$ and by
\begin{align}
q_{\text{recon}}(\vec x) = q_{\text{recon},ij}(\vec x - \vec x_{ij})  \qquad \vec x \in \mathcal C_{ij} \label{eq:recononecellonly}
\end{align}
the reconstruction on the entire grid as introduced at the end of section \ref{sec:fv}. 

In this coordinate frame, the locations of the {point values at the cell boundary} are denoted by $\vec x_m$, $m = 1, \ldots, 8$. E.g. $\vec x_4 = (\Delta x /2, 0)$.

To simplify notation it is also helpful to define a mapping to a reference cell. Consider coordinates $\vec x = (x, y)$ aligned with the cell edges and centered at a given cell. They can be transformed to reference coordinates $\vecxi = (\xi, \eta)$ of a square cell with dimensions $\Delta x=\Delta y=2$. The edges of the square are located at $\xi=-1$, $\xi=1$, $\eta=-1$ and $\eta=1$. The transformation between the coordinates is
\begin{align}
\vecxi&=
 J^{-1}
 \vec{x} &&\text{with} & 
  J^{-1} &= \begin{pmatrix}
\frac{2}{\Delta x} & 0 \\
0 & \frac{2}{\Delta y} \\
\end{pmatrix} \label{eq:jacobian}
\end{align}

$\vecxi_m$ are the locations of the {point values at the cell boundaries} in the reference cell. They correspond to $\vec x_m = J \vecxi_m$. The reconstruction $q_{\text{recon},ij}(\vec x(\vecxi))$ will be denoted by the same symbol $q_{\text{recon},ij}(\vecxi)$ whenever there is no confusion possible.

A straightforward derivation of the reconstruction polynomial is to solve the linear system that arises by inserting a general biquadratic polynomial
\begin{align}
a_{00} + a_{10} \xi + a_{20} \xi^2 + a_{01} \eta + a_{11} \xi \eta + a_{21} \xi^2 \eta + a_{02} \eta^2 + a_{12} \xi \eta^2 + a_{22} \xi^2 \eta^2 \label{eq:reconpolynomial}
\end{align} 
into the 9 equations that arise from \eqref{eq:conservation} and \eqref{eq:boundaryinterpol}. The polynomial consists of all possible combinations of $\xi$ and $\eta$ each up to second degree. It has $9$ free parameters. Hence, the solution of the linear system is unique.

This reconstruction can be obtained in a more elegant way, reminiscent of Lagrange polynomials. This representation of the interpolation will also be useful later. In every cell, one seeks a reconstruction in the form
\begin{align}
q_{\text{recon},ij}(\vecxi)=\sum_{m=1}^{9}c_{m,ij} b_{m}(\vecxi)
\label{reconSum2dRect}
\end{align}
with basis functions $b_{m}(\vecxi)$ and coefficients $c_{m,ij} \in \mathbb R^n$. The expressions for $b_m$ are the same for every cell while the coefficients $c_{m,ij}$ vary, as they depend on the average value in the cell and the {point values at the cell boundaries}.

The reconstruction can be organized nicely by choosing the basis functions at most biquadratic with
\begin{align}
 b_{m}(\vecxi) &= \begin{cases} 1 & \vecxi = \vecxi_m \\ 0 & \text{else}  \end{cases} \qquad m = 1, \ldots, 8 \label{eq:interpolprop}\\
 b_{9}(\vecxi_m) &= 0 \quad \forall m \in \{ 1, \ldots, 8\} \label{eq:interpolprop9}\\
 \frac1{|\mathcal C|}\int_\mathcal C \dd \vec x \, q_{\text{recon}}(\vecxi(\vec x)) &= q_\mathcal C \label{eq:ninthbasisfctaverage}
\end{align}
This implies that
\begin{align}
 c_{m,ij} = q_{m,ij} \qquad \forall m \in \{ 1, \ldots, 8\} \label{eq:valuec18}
\end{align}
The value of $c_{9,ij}$ is to be determined after the basis functions have been obtained explicitly.

\begin{theorem}[Interpolation basis] \label{reconRecThm}
 The following polynomials fulfill \eqref{eq:interpolprop}--\eqref{eq:interpolprop9}
\begin{align}
b_1 &= -\frac{1}{4}(\xi-1)(\eta-1)(\eta+\xi+1)\\ 
b_2 &= \frac{1}{2}(\xi-1)(\eta-1)(\xi+1)		\\ 
b_3 &= \frac{1}{4}(\xi+1)(\eta-1)(\eta-\xi+1)  \\ 
b_4 &= -\frac{1}{2}(\eta-1)(\xi+1)(\eta+1)	\\ 
b_5 &= \frac{1}{4}(\xi+1)(\eta+1)(\eta+\xi-1)  \\ 
b_6 &= -\frac{1}{2}(\xi-1)(\eta+1)(\xi+1)		\\ 
b_7 &= -\frac{1}{4}(\xi-1)(\eta+1)(\eta-\xi-1) \\ 
b_8 &= \frac{1}{2}(\eta-1)(\xi-1)(\eta+1)		\\
b_9 &= (\eta-1)(\eta+1)(\xi-1)(\xi+1)  
\end{align}
\end{theorem}
\begin{proof}
 The result is verified by explicit calculation.
\end{proof}
\emph{Note:} These products are polynomials that vanish along straight lines in $\xi$--$\eta$-plane, as shown in Figures \ref{reconEdge} and \ref{reconNode} for $b_1$ and $b_2$. 

Additionally, using \eqref{eq:valuec18} and computing \eqref{eq:ninthbasisfctaverage} one obtains
\begin{align}
 c_{9,ij} = \frac{9}{16}\Big(4 q_{ij}+  \frac{1}{3} (q_{1,ij}-4 q_{2,ij}+q_{3,ij}-4 q_{4,ij}+q_{5,ij}-4 q_{6,ij}+q_{7,ij}-4 q_{8,ij})\Big) \label{eq:value9}
\end{align}

\eqref{reconSum2dRect} becomes a biquadratic polynomial upon inserting the basis functions $b_1$, \ldots $b_9$. With the coefficients \eqref{eq:valuec18} and \eqref{eq:value9} it must be equation \eqref{eq:reconpolynomial} because the interpolation polynomial in this case is unique.

\doubleImage{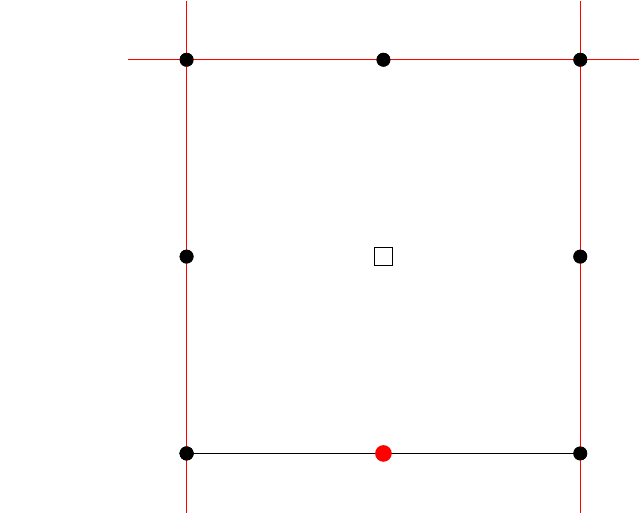}{The linear terms for the reconstruction polynomials for the location of an edge value (red point).}{reconEdge}{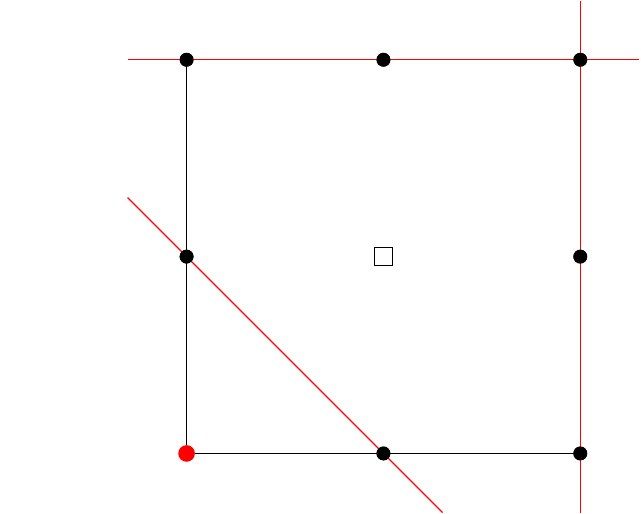}{The linear terms for the reconstruction polynomials for the location of a node value (red point).}{reconNode}

\begin{theorem}[Continuous reconstruction]
 The biparabolic reconstruction \eqref{reconSum2dRect} with basis functions from Theorem \ref{reconRecThm} and coefficients from \eqref{eq:valuec18} and \eqref{eq:value9} is continuous across cell interfaces.
\end{theorem}
\begin{proof}
From \eqref{eq:reconpolynomial} it is obvious that the reconstruction is parabolic along the coordinate axes. Consider an edge common to two cells. The two reconstructions on both sides of the edge both reduce to quadratic functions along the edge. Moreover, they agree in three points (at the two vertices and at the edge midpoint). This information specifies the quadratic function uniquely, and the two quadratic parabolas must agree.
\end{proof}

Note that the derivatives perpendicular to the common cell edge in general are different. An example of a reconstruction for an arbitrary cell is sketched in Figure \ref{reconIllustration}.

\doubleImage{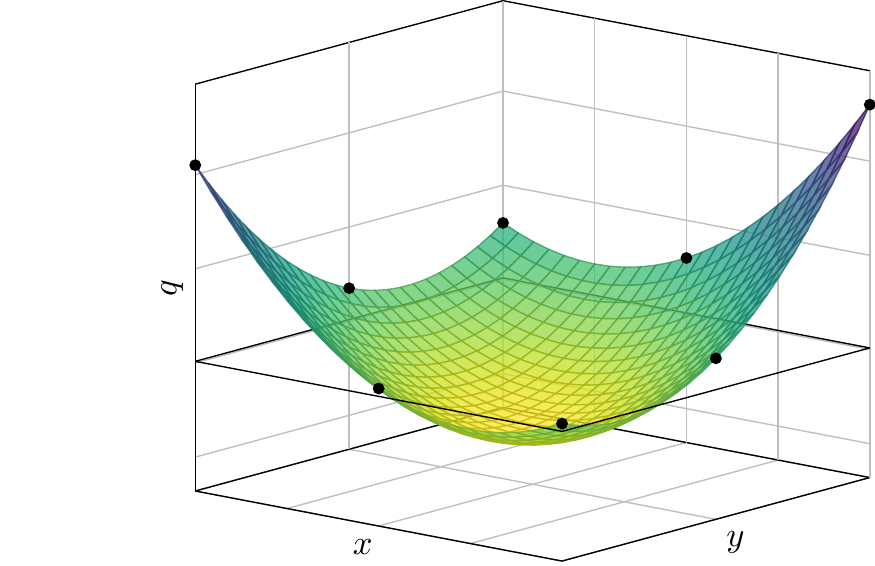}{The {point values of $q$ at cell boundaries} (black dots) and the according reconstruction of $q$.}{reconIllustration}{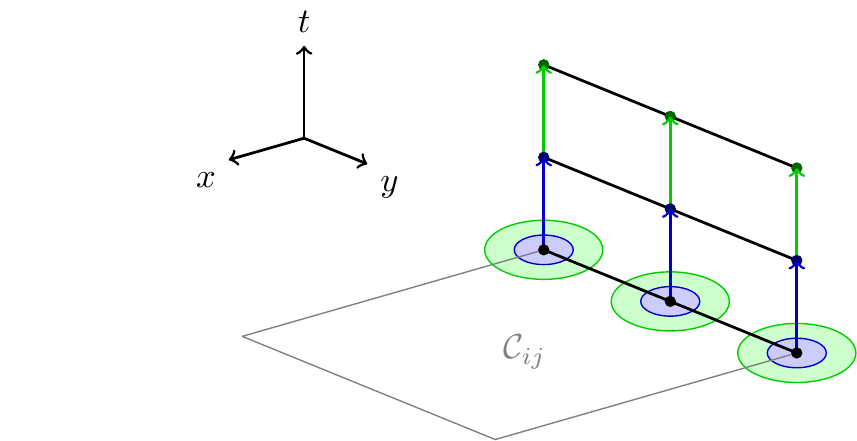}{The blue circles represent the time evolutions with the method of spherical means for the half time step; the green circles represent the time evolution with the method of spherical means for one time step. They show the domain of dependence of the solution.}{SMupdateIllustration}

\subsection{Updating the {point values at cell boundaries}}\label{chapterRectangularUpdate}

In order to advance in time the {point values located at cell boundaries}, one solves the initial value problem \eqref{eq:conslaw} with initial data
\begin{align}
 q(0, \vec x) = q_\text{recon}(\vec x)
\end{align}
at the location of the degree of freedom.

The exact solution of the acoustic equations on a rectangular mesh is given by equations \eqref{eq:solp}--\eqref{eq:solu}. The reconstruction $q_\text{recon}$ is defined piecewise according to \eqref{eq:recononecellonly}. The evolution of the {point values at the cell boundaries} therefore involves the reconstruction in several adjacent cells.

For example, the integration for the spherical mean of an edge value is split up into integrations over two hemispheres: The integration bounds for a vertical ($x = \const$) edge are $\phi\in\{\frac{\pi}{2},\frac{3\pi}{2}\}$ and $\phi\in\{0,\pi\}$ for a horizontal ($y = \const$) edge. It is convenient to define the partial spherical mean
\begin{align}
\begin{split}
	M_{\phi_1}^{\phi_2} \left[f\right](\vec{x}, r)&=\frac{1}{4 \pi}\int\displaylimits_{\phi_1}^{\phi_2} \dd \phi \int\displaylimits_0^\pi \dd \theta \, \sin \theta  f\left(\vec{x}+r\cdot \vec{n}\right)
\end{split}
\label{SM3Dpartial}
\end{align}
and $M_0^{2\pi} \equiv M$. Thus for a vertical edge, the spherical mean is computed as
\begin{align}
 M_0^{2\pi}\left[q_\text{recon}\right]\left(\vecc{x_{i+\frac12}}{y_j}, r\right) &= M_{\pi/2}^{3\pi/2} \left[q_{\text{recon},ij}\right]\left(\vecc{\Delta x / 2}{0}, r\right) \\&\phantom{mmm}+ M_{-\pi/2}^{\pi/2} \left[q_{\text{recon}, i+1,j}\right]\left(\vecc{-\Delta x/2}{0}, r\right) \label{eq:twosphmeansev}
\end{align}

Analogously the integration for the spherical mean of a node value is divided into four different integrals with integration bounds $\phi\in\{0,\frac{\pi}{2},\pi,\frac{3\pi}{2}\}$.

In a practical implementation the computational cost can be reduced by precomputing (analytically)
\begin{align*}
	M_{\phi_1}^{\phi_2} \left[x^\ell y^m \right](\vec{x}, r)
\end{align*}
for all relevant values of $\ell, m$, at all relevant locations $\vec x_m$, $m = 1, \ldots, 8$ and for all relevant angular domains. The spherical means are then obtained as linear combinations of these precomputed values. {As equation \eqref{eq:sphmeanpoly} shows, }they are polynomials in $r$, and this even allows to compute the derivatives and integrals with respect to $r$ in \eqref{eq:solp}--\eqref{eq:solu} analytically. {All the operations like spherical means and the exact evolution operator thus can be precomputed in such a way that the final formula is merely a linear combination of the values in neighbouring cells.} Thus, the evolution operator can be programmed in a way that makes its evaluation inexpensive {as compared to the evaluation of an approximate Riemann solver for a standard finite volume scheme}.

\subsection{Time integration and numerical flux}\label{chapterRectangularFlux}

In order to advance the cell average in time the fluxes through the cell boundaries are required. The strategy proposed here follows closely that of \cite{eymann13,eymann13a} on triangular grids.

\begin{figure}[ht]%
\centering
\includegraphics[width=.7\textwidth]{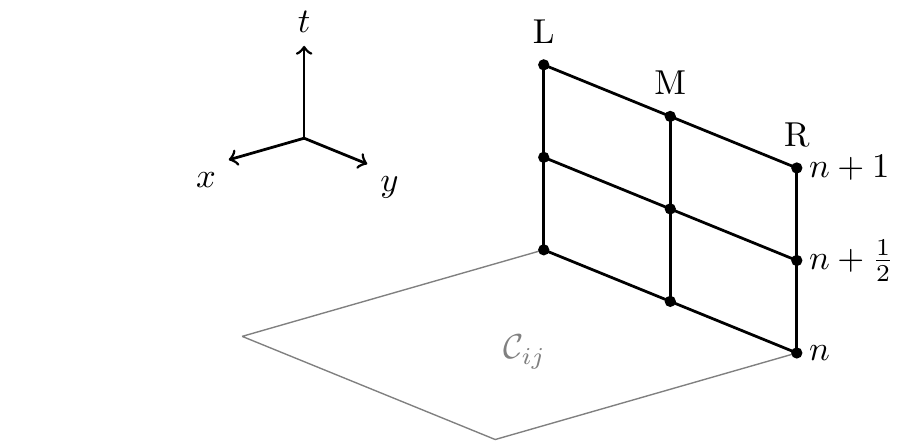}
\caption{The nomenclature for the required edge and node of quantity $q$ at different times to calculate the average flux for one side of the rectangle.}%
\label{flux2d1sRect}%
\end{figure}

Consider notation of Figure \ref{flux2d1sRect}, where an arbitrary edge is depicted and the edge value is denoted by M and the two node values by L and R.

The biparabolic reconstruction implies formally third order accuracy of the evolution of {point values along the cell boundary}. It is natural to strive for the same order of accuracy in other parts of the scheme, thus making it overall third order.

For the approximation of the average flux through a cell edge $e$
\begin{align}
 \frac{1}{|e| \Delta t} \int^{t^{n+1}}_{t^n} \!\!\!\! \dd t \int_{e} \dd x \,\vec n_e \cdot \vec f    \label{eq:exactfluxaverage}
\end{align}
two node values and one edge value are available. {As the solution is thought of as continuous across cell boundaries, these values can be immediately used in the flux function $\vec f(q)$.} Then Simpson's rule as a quadrature formula leads to a scheme of third order.

In order to achieve the same accuracy in time, it is necessary to introduce an intermediate time step $n+\frac12$. The {point values at cell boundaries} are computed at two time steps $n+\frac12$ and $n+1$. {Both go back to the initial data at time step $n$, see the green an blue regions in Figure \ref{SMupdateIllustration}.} 
{All fluxes are obtained by applying the flux function of the PDE to the point values along the boundary of the cell. Thus $f^n$ only depends on $q^n$, $f^{n+\frac12}$ only on $q^{n+\frac12}$, $f^{n+1}$ only on $q^{n+1}$, which, themselves, have been computed from the values $q^n$ during the update of the point values.}

The space-time Simpson's rule approximating \eqref{eq:exactfluxaverage} is then given by
\begin{align}
\begin{split}
 \bar f_e = \frac16 \left( \frac16 {\vec f}^n_\text L + \frac23 {\vec f}^n_\text M +  \frac16 {\vec f}^n_\text R\right )\cdot \vec n_e + \frac23 \left(\frac16 {\vec f}^{n+\frac12}_\text L + \frac23  {\vec f}^{n+\frac12}_\text M + \frac16 {\vec f}^{n+\frac12}_\text R\right )\cdot \vec n_e \\+ \frac16 \left(\frac16 {\vec f}^{n+1}_\text L +  \frac23 {\vec f}^{n+1}_\text M + \frac16 {\vec f}^{n+1}_\text R \right ) \cdot \vec n_e 
 \end{split} \label{eq:fluxbisimpson}
\end{align}

The numerical flux for a rectangular cell is visually presented in Figure \ref{flux2d1sRect}. The cell average for $t^{n+1}$ is calculated with the finite volume method \eqref{eq:fv}.

As the time integration amounts to an explicit scheme, there is a CFL-type time step restriction. All spherical means involved in the evolution of {point values at cell boundaries} have to remain inside the cell. The strongest constraint comes from the edge midpoint: the radius of the corresponding spherical mean has to be smaller than
\begin{align}
d_{\min}=\min\left(\frac{\Delta x}{2},\frac{\Delta y}{2}\right)
\end{align}

Thus the largest possible time step of the active flux method for the acoustic equations on a rectangular mesh is
\begin{align}
\Delta t_{\max}=\min \left(\frac{\Delta x}{2c},\frac{\Delta y}{2c}\right)
\label{tMaxCCRect}
\end{align}
{Experimentally, the active flux method was observed to be stable with time steps very close to this bound, suggesting that the bound is sharp.}

{As can be seen from Figure \ref{stationarityMesh}, the update of a cell average thus involves the reconstruction inside the cell and the reconstruction in all its 8 neighbours.}
The numerical fluxes and the cell average for an arbitrary cell are sketched in Figure \ref{cellAverageUpdateIllustration}.

\doubleImage{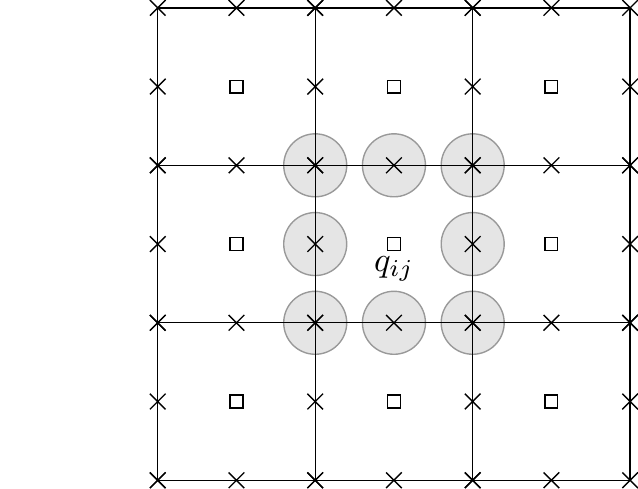}{The spheres for the spherical means, which are necessary to calculate the time update of the $ij$-th cell average ${q}_{ij}$.}{stationarityMesh}{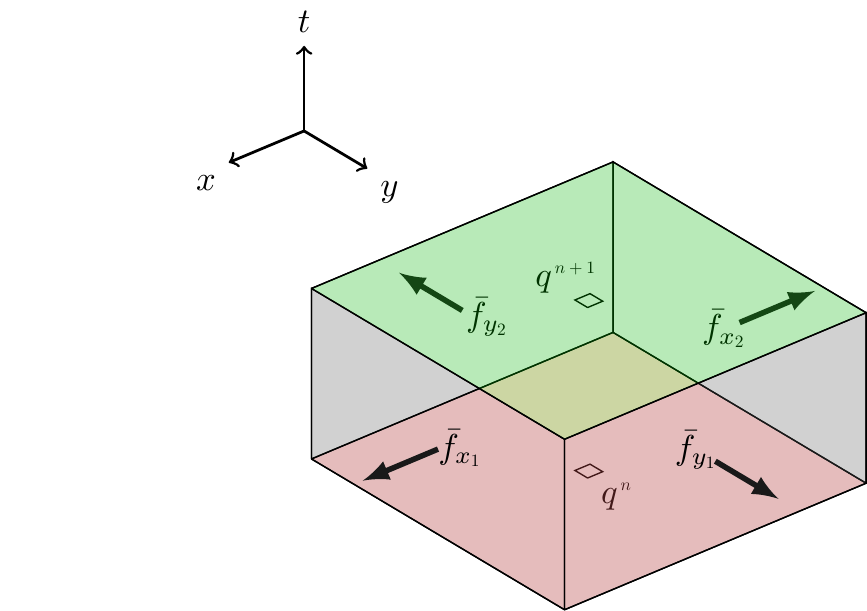}{The cell average $q^{\protect\scalebox{0.5}{$\,n$}}$ and the numerical fluxes $\bar{f}_{x_{\protect\scalebox{0.5}{$1$}}}$, $\bar{f}_{x_{\protect\scalebox{0.5}{$2$}}}$, $\bar{f}_{y_{\protect\scalebox{0.5}{$1$}}}$, and $\bar{f}_{y_{\protect\scalebox{0.5}{$2$}}}$ required to calculate the new cell average $q^{\protect\scalebox{0.5}{$\,n+1$}}$.}{cellAverageUpdateIllustration}

\section{Numerical results}\label{chapterRectangularNumerical}

In this section numerical results of the active flux scheme for the acoustic equations are presented. {Given analytical expressions for the initial data, the {point values at cell boundaries} are initialized by evaluating these expressions at the corresponding locations. The initial cell average is computed with a quadrature formula following Simpson's rule.}

\subsection{Stationary vortex}

The first setup is the stationary vortex ($r = \sqrt{x^2 + y^2}$)
\begin{align}
	\vec{v}(r) &= \vec{n}_{\phi}\cdot
		\begin{cases}
			5r & \text{for } 0 \leq r \leq 0.2  \\
			2-5r  & \text{for } 0.2 < r \leq 0.4 \\
			0 & \text{for } r > 0.4 
		\end{cases} &
			p(r) &= p_0
	\label{accousticGreshoRec}
\end{align} 
with $\vec n_{\phi}=(- \sin(\phi), \cos(\phi))^T$. The simulation has been performed on a $50\times50$ Cartesian mesh with $\Delta x = \Delta y = 0.03$. {Here, $p_0 = 0$ and the CFL number is $0.45$.}

\begin{figure}[!bt]%
\begin{minipage}{0.48\textwidth}
\includegraphics[width=\columnwidth]{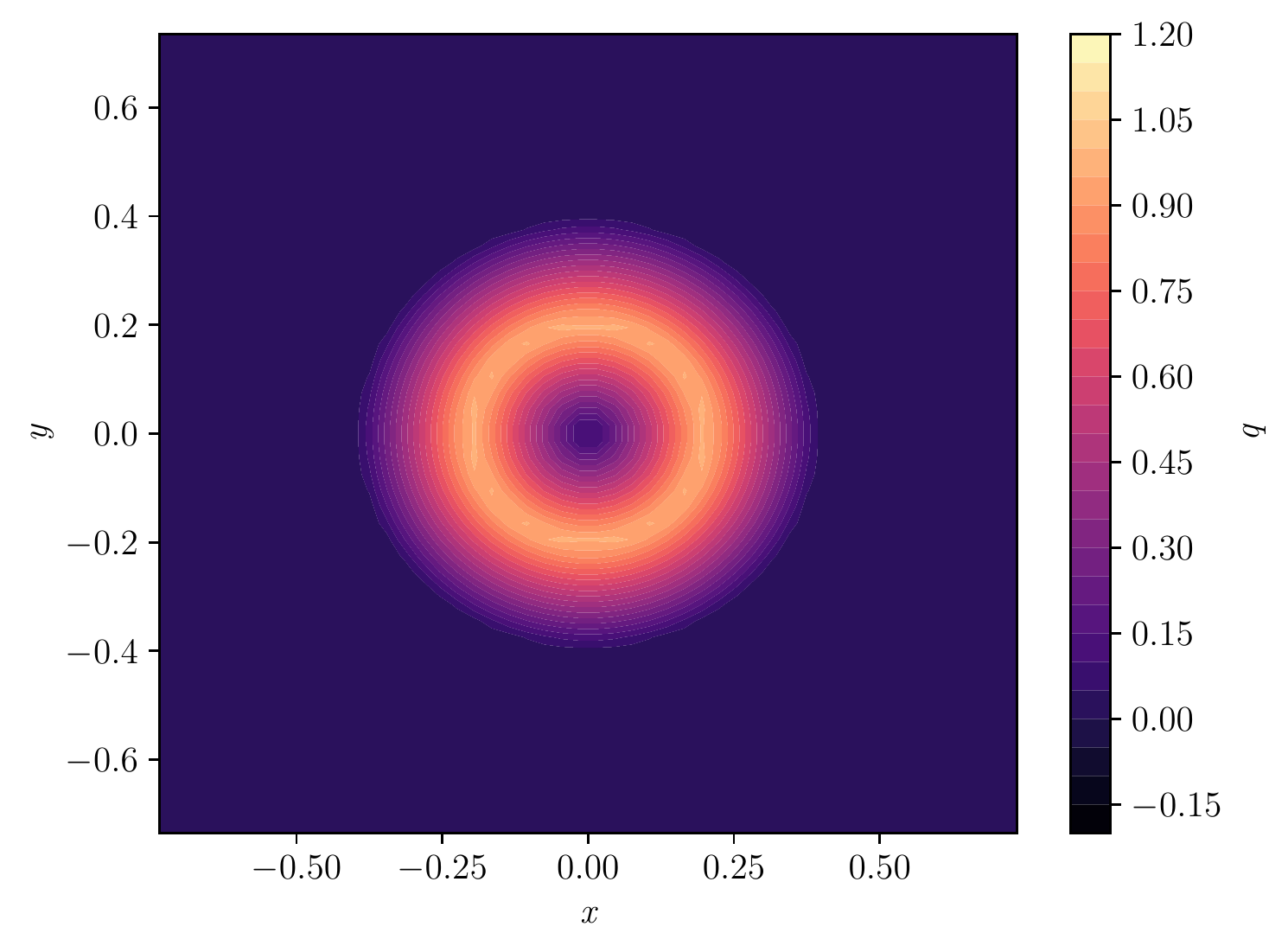}%
\caption{The vortex setup for the acoustic equations at $t=0$. Contour plots of the cell averages of $|\vec v|$ are shown.}
\label{greshoAccRec08-initial}
\end{minipage}
\hfill
\begin{minipage}{0.48\textwidth}
\includegraphics[width=\columnwidth]{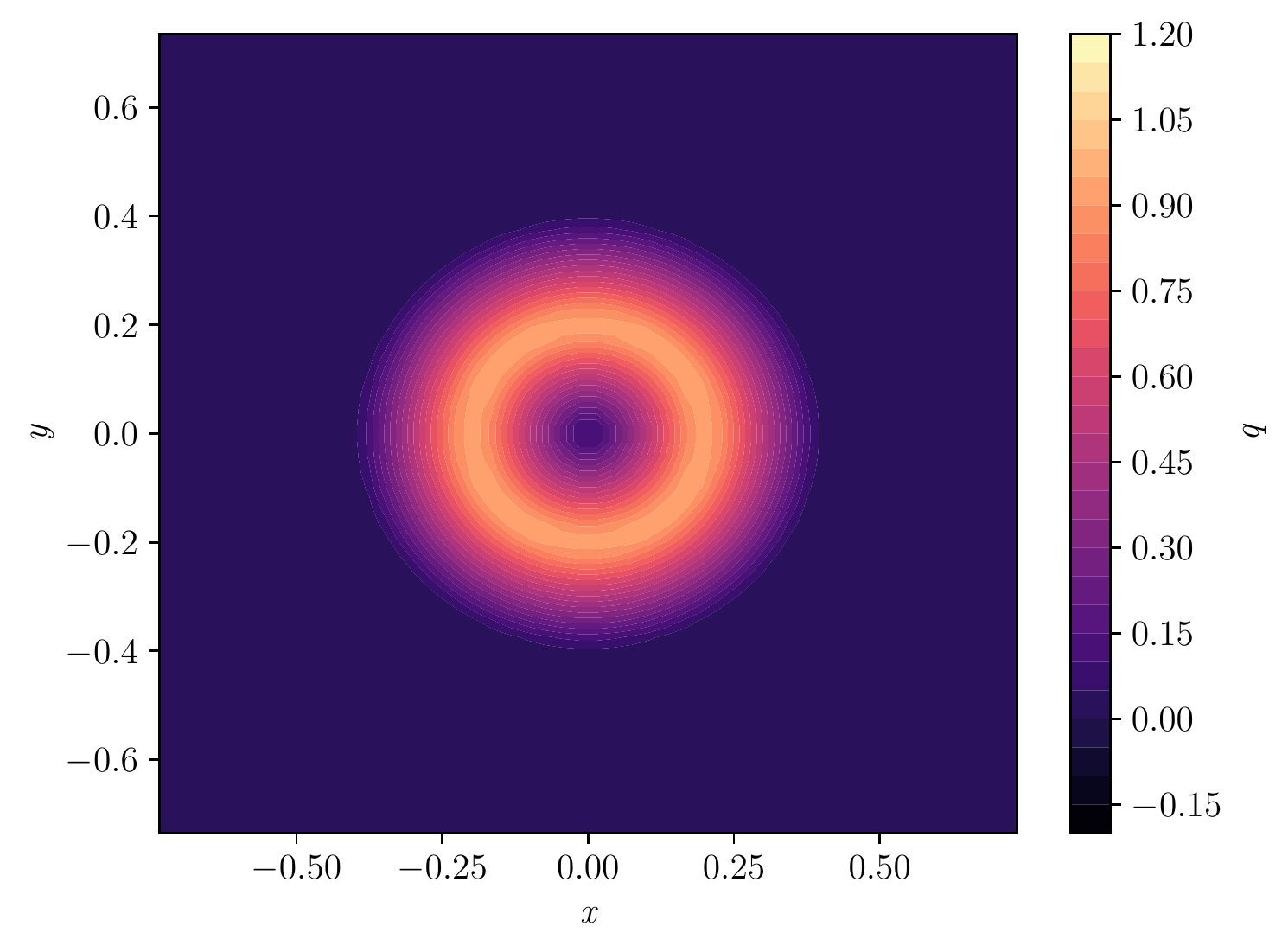}%
\caption{The vortex setup for the acoustic equations at $t=100$. Contour plots of the cell averages of $|\vec v|$ are shown.}
\label{greshoAccRec08}
\end{minipage}
\end{figure}
Contour plots of the cell averages of the absolute value of the velocity $\left|\vec{v}\right|=\sqrt{u^2+v^2}$ for this setup are shown in Figures\ref{greshoAccRec08-initial}--\ref{greshoAccRec08} for the times $t=0$ and $t=100$.  

\begin{figure}[h]%
\includegraphics[width=0.7\columnwidth]{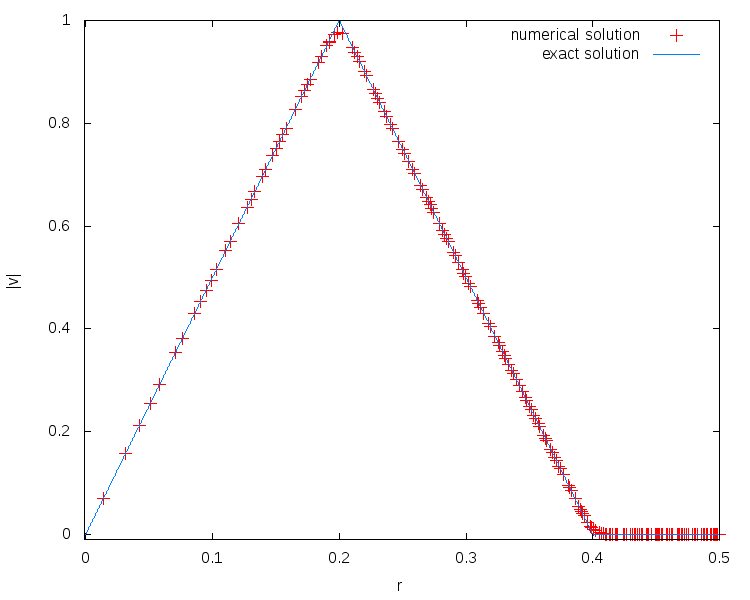}%
\caption{Radial plot of the cell averages at time $t=100$ and the initial data (solid line) of the absolute value of the velocity $\left|\vec{v}\right|$. {The depicted numerical solution represents the discrete stationary state which is maintained to machine precision (as shown in section \ref{statPresAF} by applying the theory from \cite{barsukow17a}).}}
\label{radialPlotsAccRec}%
\end{figure}

Additionally, the radial plot of the cell averages of the absolute value of the velocity $\left|\vec{u}\right|$ at $t=100$ and the exact solution are shown in Figure \ref{radialPlotsAccRec}. The radial plots in Figure \ref{radialPlotsAccRec} suggest that the presented method is able to preserve stationary states to very high accuracy. This is not just due to the high order of the scheme. As is shown in the next section, the active flux method on rectangular grids is \emph{stationarity preserving}. This is a property introduced in \cite{barsukow17a} which means that the stationary states of the scheme discretize all the analytic stationary states.

\begin{figure}[h]
 \includegraphics[width=0.49\textwidth]{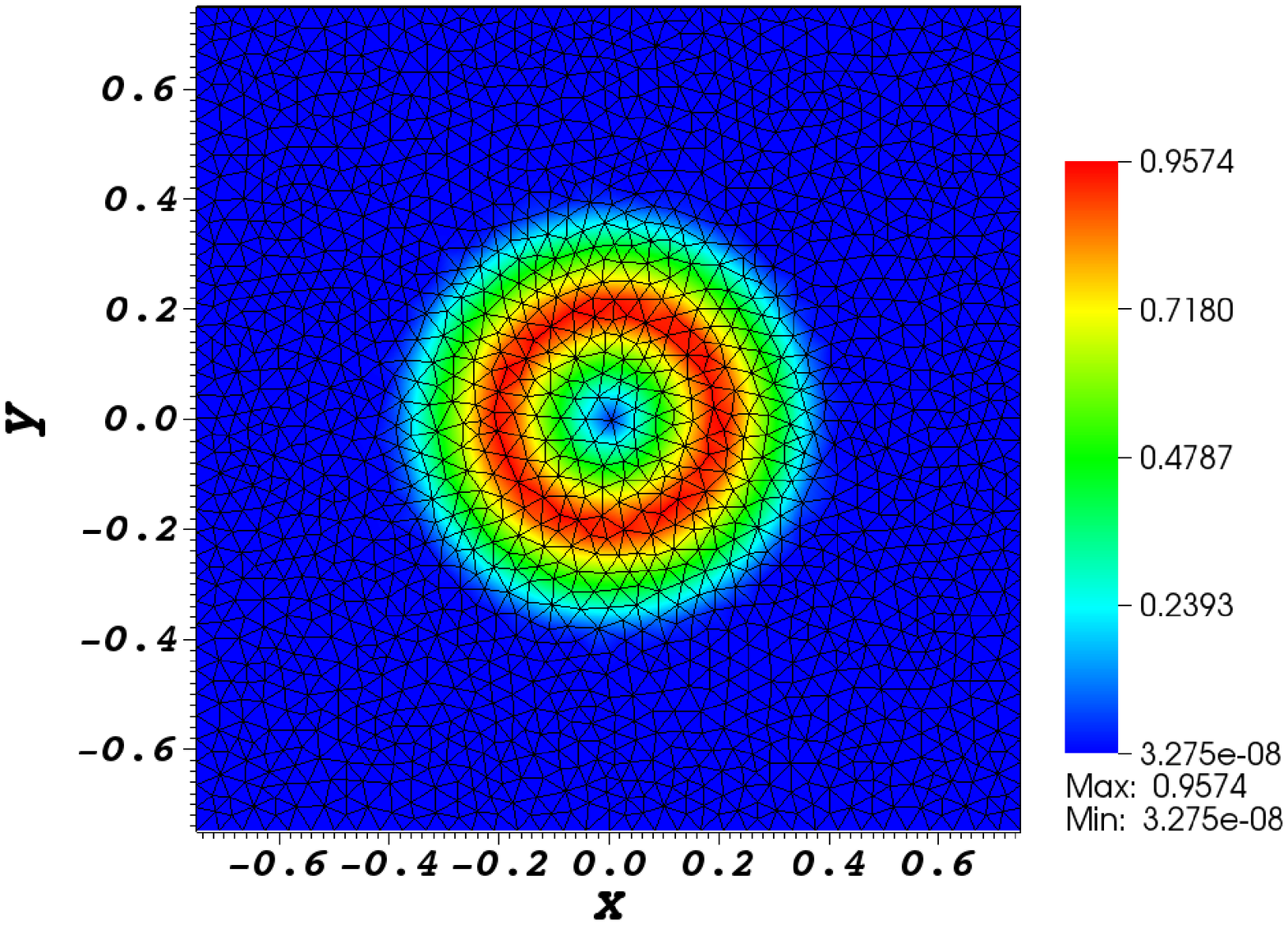}\hfill\includegraphics[width=0.49\textwidth]{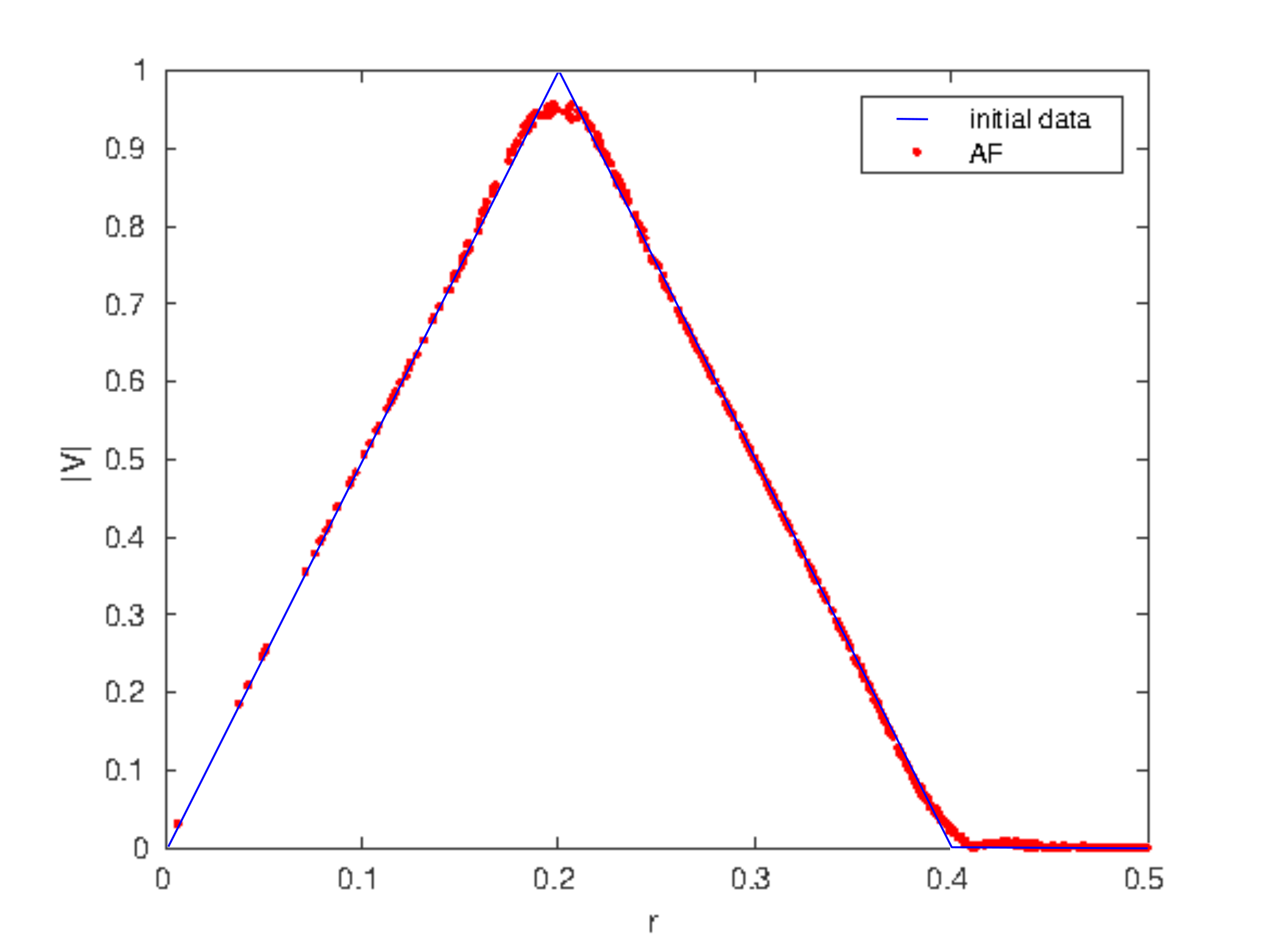}
 \caption{{The vortex setup for the acoustic equations solved using the active flux method from \cite{eymann13} on a 
 triangular grid (see section \ref{sec:standardaf}). One observes a behaviour comparable to the one on Cartesian grid.}}
 \label{fig:unstructuredgresho}
\end{figure}

{Additionally, Figure \ref{fig:unstructuredgresho} shows the same setup on a grid with a similar number of triangular elements. Average values for every element are shown.}

\subsection{Acoustic waves}

\begin{figure}[h]
 \centering
 \includegraphics[width=0.48\textwidth]{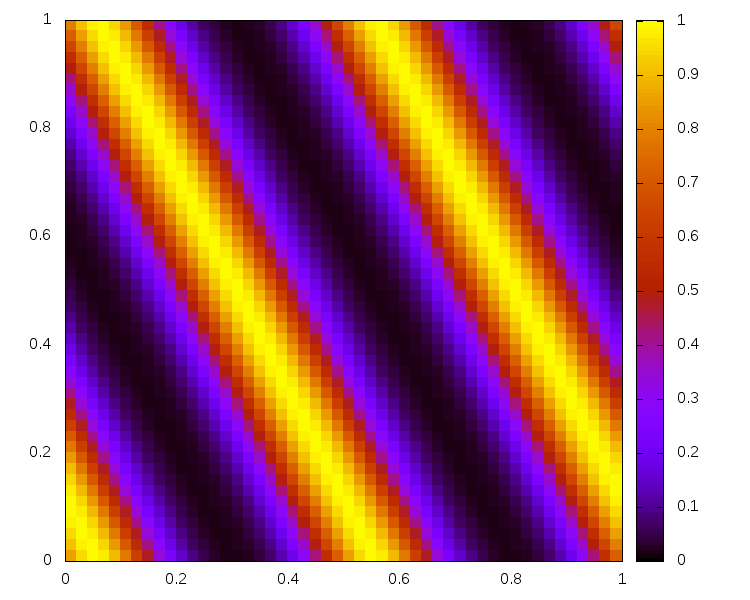}\hfill\includegraphics[width=0.48\textwidth]{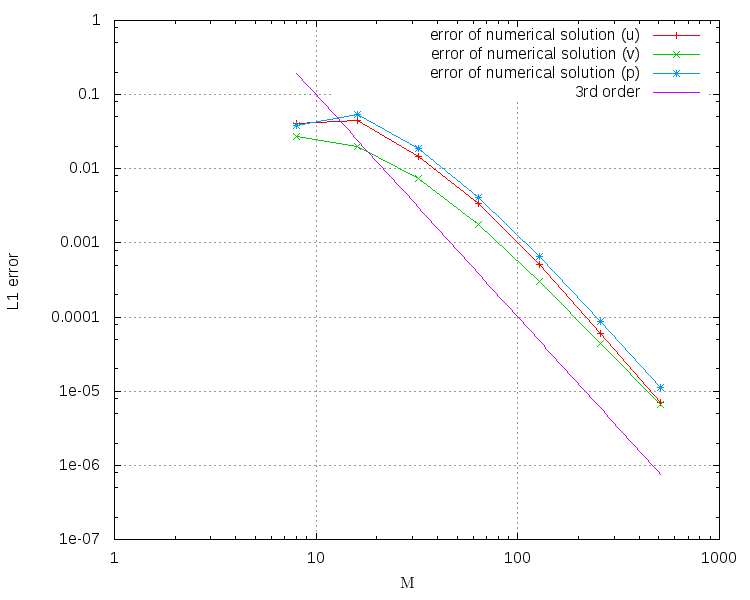}\\
 \includegraphics[width=0.48\textwidth]{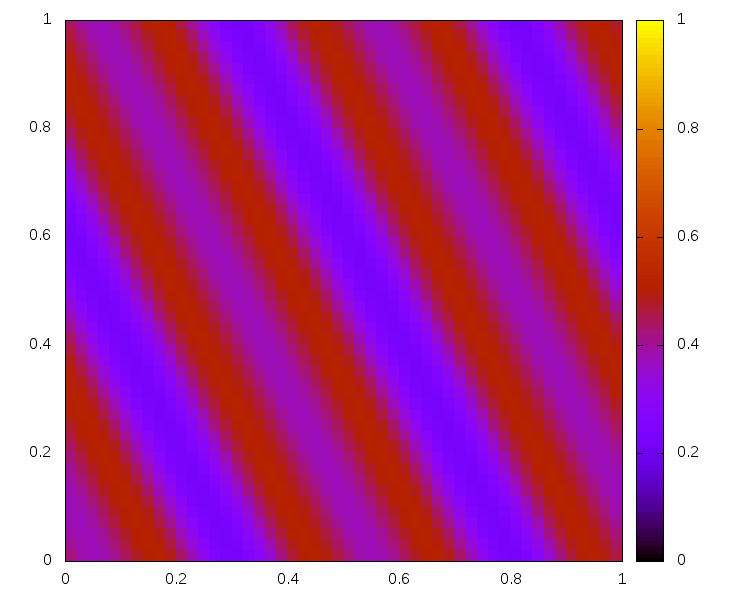}\hfill\includegraphics[width=0.48\textwidth]{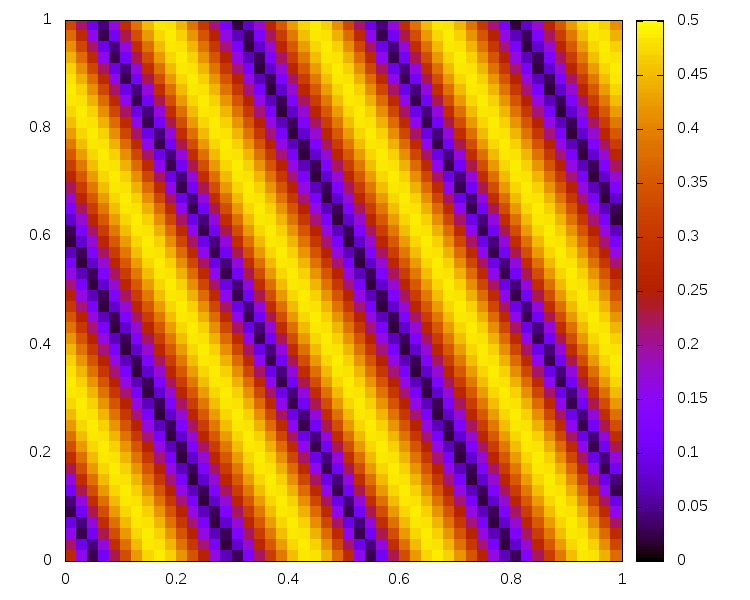}
 \caption{Acoustic waves propagating obliquely to the grid directions. \textit{Upper row, left}: Pressure at $t=0$ on a $50\times50$ grid. \textit{Upper row, right}: $L^1$ error of the solution at vertical edge midpoints at $ct = 0.1$ for different $M \times M$ grids. A third order convergence is observed. \textit{Lower row}: Pressure (\textit{left}) and the absolute value $|\vec v|$ (\textit{right}) of the velocity at $ct = 0.1$ on a $50\times50$ grid.}
 \label{fig:gaussianwave}
\end{figure}

A second test case is that of smooth acoustic waves propagating obliquely to the grid directions:
\begin{align}
	\vec{v} &= 0 &
	p(\vec x) &= \sum_{i=-2}^2 \exp\left( - \frac{(\vec x \cdot \vec b-i\delta)^2}{w^2} \right )
	\label{gaussianwavesinit}\\
	&& \vec b &= \vecc{\,\,\cos \alpha\,\,}{\sin \alpha}
\end{align} 
Here $\delta = 0.1$, $w = 0.5 \cos \alpha$ and $\alpha = \arctan \frac12$ is used. The initial setup on a $50 \times 50$ grid covering $[0,1]^2$ with periodic boundaries is shown in Figure \ref{fig:gaussianwave} (upper left). This problem can be transformed into a one-dimensional setup along $\vec b$ which allows to derive an exact solution by the method of characteristics. Figure \ref{fig:gaussianwave} (upper right) shows the error of the numerical solution at the midpoints of the vertical edges on $M \times M$ grids for different values of $M$. For this instationary setup a CFL number of 0.45 has been used. One observes a third order convergence of the method.

\subsection{Spherical shock tube}

\begin{figure}[h]
 \centering
 \includegraphics[width=0.45\textwidth]{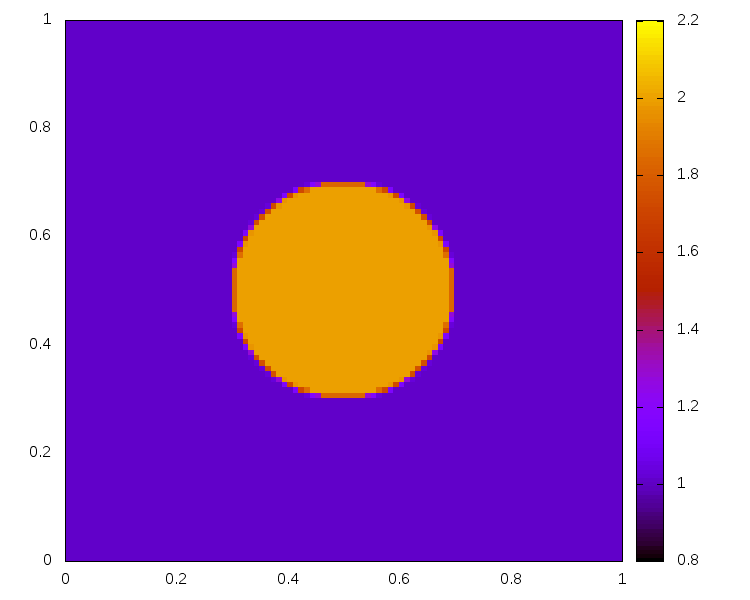}\hfill\includegraphics[width=0.45\textwidth]{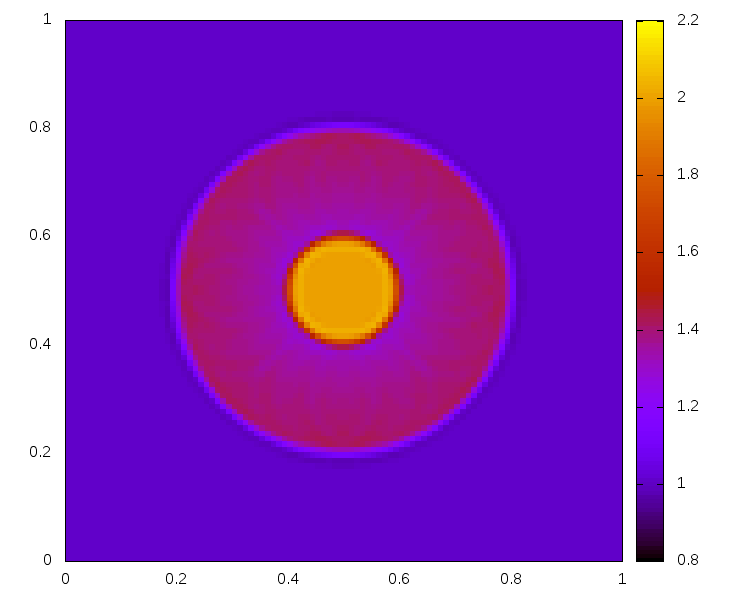}\\
 \includegraphics[width=0.45\textwidth]{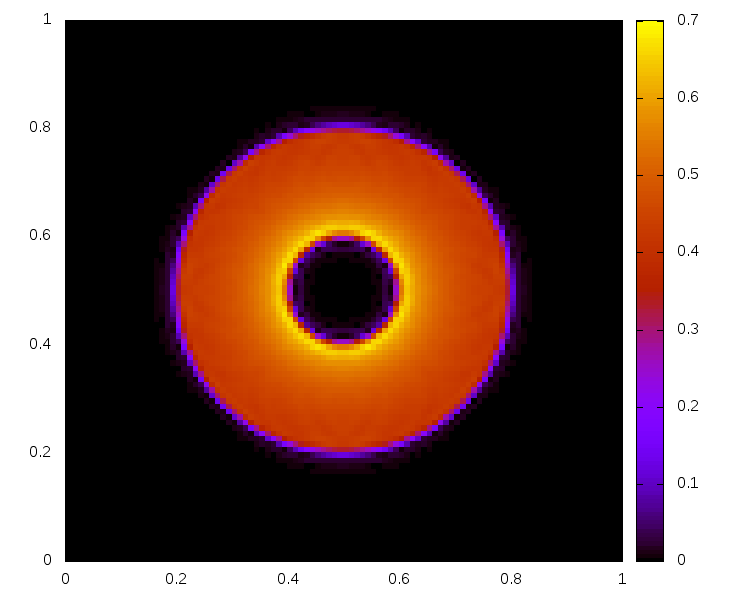}\hfill\includegraphics[width=0.45\textwidth]{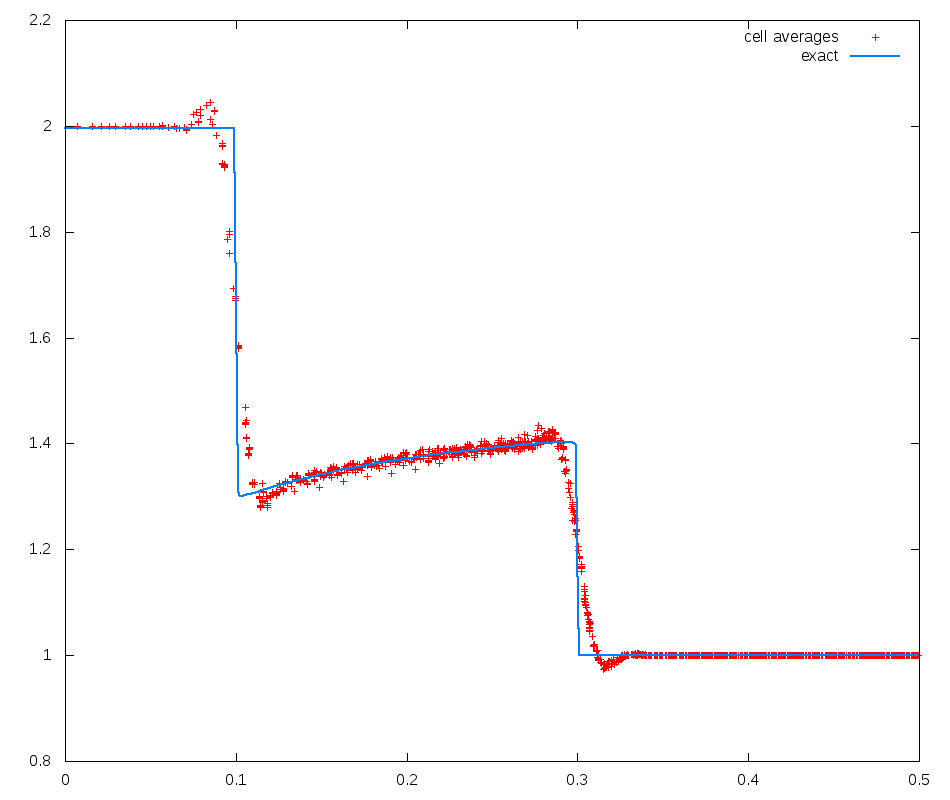}\\
 \includegraphics[width=0.45\textwidth]{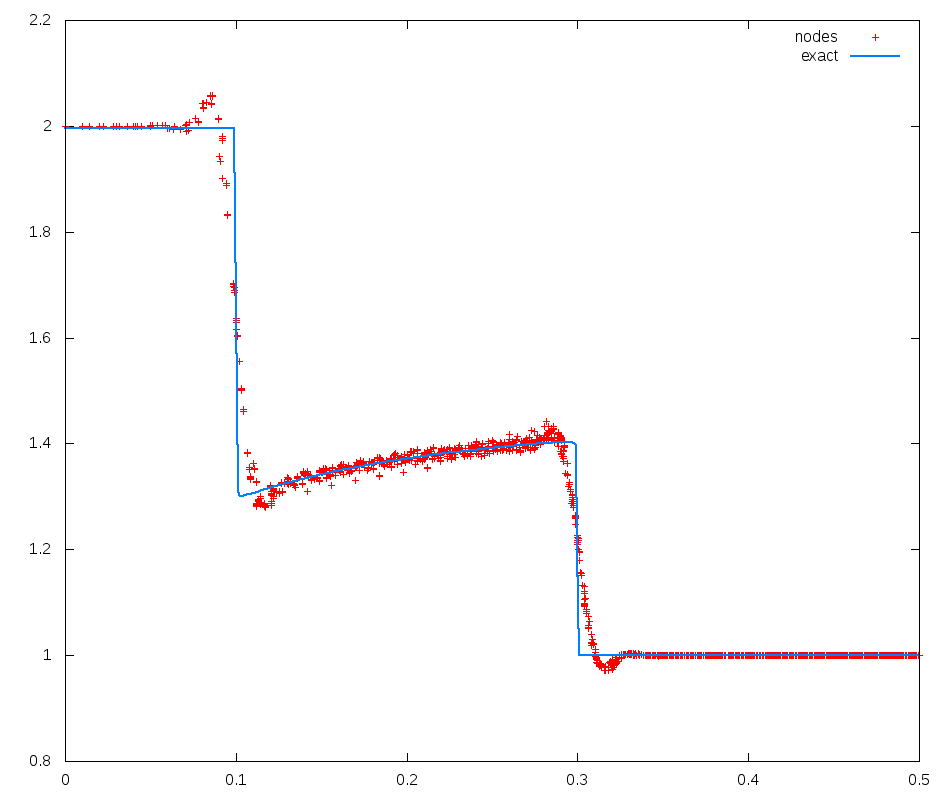}\hfill\includegraphics[width=0.45\textwidth]{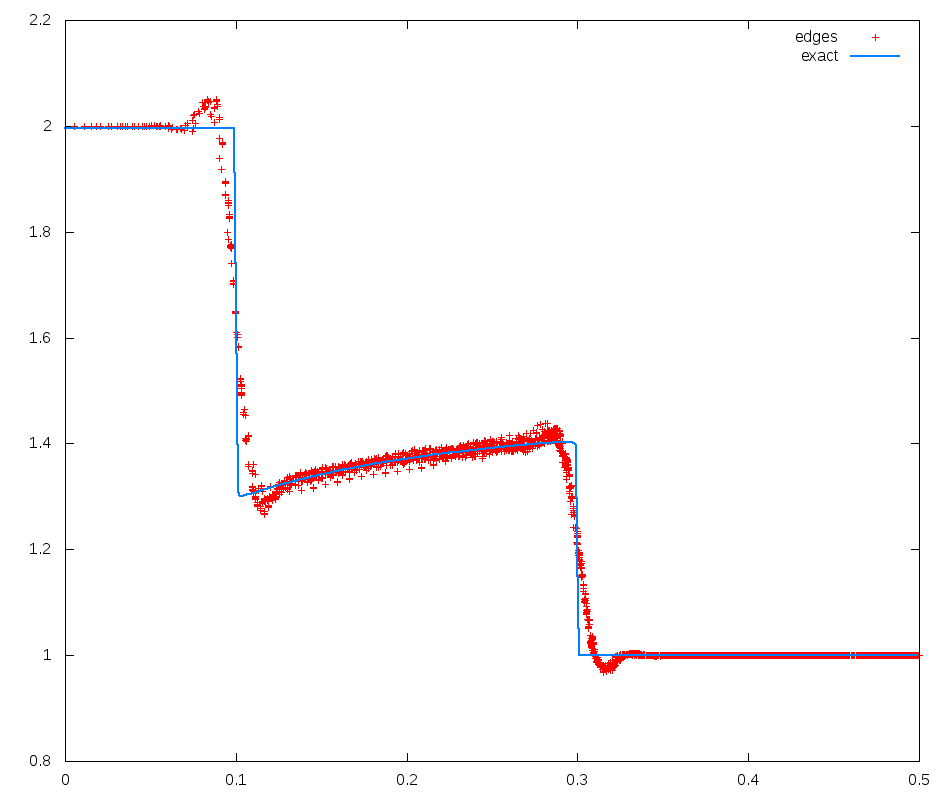}
 \caption{Radial Riemann problem solved with the active flux method on a $100\times100$ grid. The initial data involve a symmetric discontinuity in the pressure $p$ and $\vec v = 0$ everywhere. \textit{Upper row}: Initial data in the pressure (\textit{left}) and the pressure at $ct=0.1$ (\textit{right}). \textit{Center row}: Absolute value $|\vec v|$ of the velocity at time $ct=0.1$ (\textit{left}) and a scatter plot of the pressure (cell averages) at $ct=0.1$ as a function of the radius (\textit{right}). {The reference solution (solid line) has been computed by numerical quadrature of \eqref{eq:solp}.} \textit{Lower row}: Analogous scatter plots showing the values at the nodes (\textit{left}) and the edges (\textit{right}).}
 \label{fig:shock}
\end{figure}

The last test case is that of a radial shock tube. The initial data are chosen
\begin{align}
	\vec{v} &= 0 &
	p(r) &= 
		\begin{cases}
			2 & \text{for } r \leq 0.2  \\
			1 & \text{else}
		\end{cases}
	\label{radialShock}
\end{align} 
The simulation is performed on a $100\times100$ grid. The results are shown in Figure \ref{fig:shock}. One observes that despite a continuous reconstruction the method is able to compute discontinuous solutions. The radial scatter plot demonstrates the symmetry of the solution. The over- and undershoots are due to the high order of the scheme and no limiting being employed.

\section{Low Mach number compliance and vorticity preservation}\label{statPresAF}

Low Mach number compliance of a numerical scheme for linear acoustics has been found to be related to further concepts in \cite{barsukow17a}. In particular, it has been shown that \emph{stationarity preserving} schemes for the acoustic equations are just those that comply with the limit of low Mach number. 

{Preserving stationary states is a challenge to numerical methods. Whereas in one spatial dimension the stationary states of \eqref{eq:conslaw} are characterized by just an ODE, the multi-dimensional stationary states themselves fulfill a PDE. In particular, stationary states of \eqref{acousticEq} are characterized by a constant pressure and a divergencefree velocity. It is impossible to require a numerical method to keep stationary any sampling a divergencefree vector field. Thus, when treating stationary states numerically in multiple spatial dimensions it is necessary to content oneself with less.

\cite{barsukow17a} made the following observation: even though the stationary states of \eqref{acousticEq} are all divergencefree vector fields, many schemes only keep trivial states stationary, such as the constant state $q = \const$. The stationary states of such numerical schemes thus are not a discretization of \emph{all} stationary states of the PDE. In \cite{barsukow17a} therefore, \emph{stationarity preserving schemes} were introduced. They possess discrete stationary states that discretize all the analytic stationary states. In view of what has been said above, such schemes of course do not keep stationary \emph{any} sampling of the analytic stationary states. However, for example in the case of linear acoustics, they keep stationary a particular \emph{discretization} of the divergence. Initial data, sampled from a divergencefree vector field, evolve in time towards one of the discrete stationary states, which then persists for an arbitrarily long time without being diffused away.} For a thorough presentation of these concepts the reader is referred to \cite{barsukow17a,barsukow17lilleproceeding,barsukow18thesis}.

It turns out that for the acoustic equations stationarity preservation is equivalent to the scheme being \emph{vorticity preserving}, i.e. possessing a discrete vorticity that remains stationary. Therefore a numerical scheme that is not vorticity preserving also fails to discretize all the stationary states. It introduces so much diffusion that all but the most trivial stationary states are decaying in time. Numerical results of section \ref{chapterRectangularNumerical} suggest that this is not the case for the active flux scheme on Cartesian grids.

In this chapter, it is shown that the active flux method on a rectangular mesh is stationarity preserving when used for numerically solving the two-dimensional acoustic equations. The framework presented in \cite{barsukow17a} has to be generalized in order to take into account the {point values at cell boundaries} that are present in the active flux scheme. 
The proof involves applying the discrete Fourier transform to the numerical scheme, which is discussed in section \ref{sec:fouriergeneral}. Stationarity preservation of the active flux scheme is proven in Corollary \ref{thm:statpreservAF} in section \ref{ssec:statpresrecon}. The discrete vorticity that is kept stationary by the scheme is discussed there as well.

\subsection{The active flux scheme in Fourier space} \label{sec:fouriergeneral}

Stationarity preservation can be most efficiently studied in Fourier space. A prerequisite is that the scheme under consideration is linear. For the active flux scheme this is established in the following

\begin{theorem}
 The active flux scheme for the acoustic equations on Cartesian grids is linear in the discrete degrees of freedom.
\end{theorem}
\begin{proof} 
 The reconstruction \ref{reconSum2dRect} is linear in the cell average and the {point values at cell boundaries}. The spherical mean (equation \ref{SM3D}) is a linear functional. Thus the update of the {point values at cell boundaries} from time step $n$ to the time step $n+\frac12$ or $n+1$ is a linear function of the quantities at time step $n$.  The numerical flux obtained by the Simpson's rule is a linear functional of the involved {point values at cell boundaries}. Finally, the finite volume method used to update the cell average is linear in the fluxes. 
\end{proof}

{As the acoustic equations are a linear problem, and thus the active flux scheme is linear as well}, it is possible to consider it in Fourier space. This is analogous to the procedure in \cite{barsukow17a}, and switching to Fourier space allows to prove stationarity preservation easily. In the following, the nomenclature from Figure \ref{stationarityGrids} in section \ref{chapterboundarydeg} is used to differentiate between the different kinds of degrees of freedom. The additional {point values at cell boundaries} that are present in the active flux scheme can be easily incorporated in the framework. The discrete Fourier transform with respect to the spatial variable can be applied to each of the lattices $q_{ij}$, $q\no_{ij}$, $q\ev_{ij}$, $q\eh_{ij}$ defined in section \ref{chapterboundarydeg}: Any function $q_{ij}(t)$ is expressed as
\begin{align}
q_{ij}(t)=\hat{q}(t)\exp(\ii  \Delta x k_x i+\ii  \Delta y k_y j)
\label{fourierAnsatz}
\end{align}
with the wave vector $\vec{k}=(k_x,k_y)$ and the imaginary unit $\ii $. The hat in $\hat{q}$ denotes the discrete spatial Fourier transform of $q$. The overall prefactor $\exp(\ii  \Delta x k_x i+\ii  \Delta y k_y j)$ will appear in all terms.

For linear acoustics solved with the active flux method, the Fourier transforms shall be ordered as follows (see also Figure \ref{stationarityGrids})
\begin{align}
\hat{Q} :=\left( \hat{p}\avg,\hat{u}\avg,\hat{v}\avg,\hat{p}\eh,\hat{u}\eh,\hat{v}\eh,\hat{p}\ev,\hat{u}\ev,\hat{v}\ev,\hat{p}\no,\hat{u}\no,\hat{v}\no\right) \in \mathbb{R}^{12}
\label{numericalQuantities12}
\end{align}

The translation operators that convey shifts by one cell in $x$-direction and $y$-direction are defined as
\begin{align}
t_x=\exp(\ii  \Delta x k_x ) \quad \text{and} \quad t_y=\exp(\ii  \Delta y k_y)
\label{TranslationOperator}
\end{align}
such that Equation \ref{fourierAnsatz} becomes
\begin{align}
q_{ij}^{\,n}=\hat{q}^{\,n} t_x^i t_y^j
\label{fourierAnsatzOperators}
\end{align}

The Fourier transform of any two-dimensional linear finite difference formula corresponds to a Laurent polynomial in $t_x$ and $t_y$. For example:
\begin{align}
 {q}_{i+1,j} &= \hat{{q}} \, t_x \cdot t_x^it_y^j & {q}_{i,j+1} &= \hat{{q}} \, t_y \cdot t_x^it_y^j
\end{align}
A more complicated finite difference formula can be, for example, expressed as
\begin{align} 
 q_{i-1} + 4 q_i + q_{i+1} = \hat q \, \left( t_x^{-1} + 4 + t_x  \right) \cdot t_x^it_y^j = \hat q \frac{1 + 4 t_x + t_x^2}{t_x} \cdot t_x^it_y^j
\end{align}

Also the reconstruction $q_{\text{recon}, ij}(\vec x)$ in cell $\mathcal C_{ij}$ is a linear function of the degrees of freedom in $\mathcal C_{ij}$ which allows to apply the discrete Fourier transform to it:
\begin{align}
 q_{\text{recon}, ij}(\vec x) = \hat q_{\text{recon}}(\vec x) \cdot t_x^it_y^j
\end{align}
Note that the spatial variable $\vec x$ is unaffected, because the origin $\vec x = 0$ has been chosen to be the cell midpoint of the corresponding cell, and the discrete Fourier transform only acts on $i,j$. Then upon the discrete Fourier transform one can rewrite the reconstruction in any other cell $\mathcal C_{i+s, j+p}$ as
\begin{align}
 q_{\text{recon}, i+s,j+p}(\vec x) = t_x^s t_y^p \hat q_{\text{recon}}(\vec x) \cdot t_x^it_y^j
\end{align}

On the other hand, the evolution operator \eqref{eq:solp}--\eqref{eq:solu} only acts on $\vec x$. It is thus possible to construct a $12 \times 12$ matrix $\mathcal A $ which describes the evolution of the Fourier modes $\hat Q$ introduced in \eqref{numericalQuantities12}:
\begin{align}
 \hat Q(\Delta t) = \mathcal A(\Delta t; t_x, t_y) \hat Q(0)
\end{align}

\begin{theorem} \label{thm:stationarityevolutionmatrix}
 The active flux scheme for linear acoustics in two spatial dimensions is stationarity preserving iff $\det \Big(\mathcal A(\Delta t; t_x, t_y) - \id_{12\times12} \Big ) = 0$ $\forall \Delta t, t_x, t_y$.
\end{theorem}
\begin{proof}
 This is Theorem 2.11 from \cite{barsukow17a}.
\end{proof}

The matrix $\mathcal A$ can be constructed explicitly. However, this involves computations of considerable length. Therefore in the following section a simpler strategy is applied. Its main ingredient is the observation that the active flux scheme uses the \emph{exact} evolution operator for the {point values at cell boundaries}.

\subsection{Stationarity preserving reconstruction} \label{ssec:statpresrecon}

Clearly, at continuous level, if the initial data for \eqref{acousticEq} fulfill
\begin{align}
 \nabla p_0 &= 0 & \div \vec v_0 &= 0
\end{align}
then they remain stationary. (This statement needs to be understood in the sense of distributions, if necessary -- see \cite{barsukow17} for more details.) Therefore one path towards understanding the stationary states of the active flux scheme is to study {under which conditions} the reconstruction fulfills
\begin{align}
 \nabla p_\text{recon} &= 0 & \div \vec v_\text{recon} &= 0 \label{eq:stationaryrecon}
\end{align}

Assume this to be the case for some choice of the discrete values on the grid. Then at least the {point values at cell boundaries} will remain precisely stationary over one time step. This shall be studied now in more detail; the question of whether the cell averages remain stationary is postponed and taken up at the end.

Theorem \ref{thm:stationarityevolutionmatrix} states that a scheme is stationarity preserving (and thus vorticity preserving) if there is a zero eigenvalue. The corresponding eigenvector can be seen as a non-zero choice of 
\begin{align}
 \hat Q = \left( \hat{p}\avg,\hat{u}\avg,\hat{v}\avg,\hat{p}\eh,\hat{u}\eh,\hat{v}\eh,\hat{p}\ev,\hat{u}\ev,\hat{v}\ev,\hat{p}\no,\hat{u}\no,\hat{v}\no\right)
\end{align}
that implies that all the degrees of freedom remain stationary. Therefore

\begin{theorem} \label{thm:stationarityrecon}
The active flux scheme for linear acoustics is stationarity preserving, if there is a non-zero choice of $\hat Q$ that implies \eqref{eq:stationaryrecon}.
\end{theorem}
\emph{Note:} As such a non-zero $\hat Q$ spans an eigenspace, any multiple of it implies \eqref{eq:stationaryrecon} as well.

Observe the tremendous simplification from Theorem \ref{thm:stationarityevolutionmatrix} to Theorem \ref{thm:stationarityrecon}. Whereas showing $\det \Big(\mathcal A(t_x, t_y) - \id_{12\times12} \Big ) = 0$ requires writing down explicitly all the spherical means (like in \eqref{eq:twosphmeansev}) and evolution operators, the statement of Theorem \ref{thm:stationarityrecon} refers to a property of the initial data only. This is only possible because the evolution operators are exact and one thus can express in simple words which data they keep stationary.

\begin{theorem} \label{thm:stationaryeigenvector}
 If $\hat Q$ is such that
 \begin{align}
   \hat{p}\avg &= \hat{p}\eh = \hat{p}\ev = \hat{p}\no = 0
  \end{align}
  \begin{align}
   \hat{u}\avg &= -\frac23 \frac{1 + 4 t_x + t_x^2}{t_x} \cdot \frac{(t_y-1)(t_y+1)}{\Delta y t_y} &
   \hat{v}\avg &= \frac23 \frac{1 + 4 t_y + t_y^2}{t_y}  \cdot \frac{(t_x-1)(t_x+1)}{\Delta x t_x} \label{eq:statinitialdatafourier1}\\
   \hat{u}\eh &= - \frac{1 + 6 t_x + t_x^2}{t_x} \cdot  \frac{t_y-1}{\Delta y}&
   \hat{v}\eh &= 2 \frac{(t_x-1) (t_x+1)}{\Delta x t_x} (t_y+1) \label{eq:statinitialdatafourier2}\\
   \hat{u}\ev &= -2 (t_x+1) \frac{(t_y-1)(t_y+1)}{\Delta y t_y} &
   \hat{v}\ev &= \frac{t_x-1}{\Delta x}  \cdot \frac{1 + 6 t_y + t_y^2}{t_y} \label{eq:statinitialdatafourier3}\\
   \hat{u}\no &= - 4 (t_x+1)\frac{t_y-1}{\Delta y}&
   \hat{v}\no &= 4 \frac{t_x-1}{\Delta x} (t_y+1) \label{eq:statinitialdatafourier4}
 \end{align}
 then $p_\text{recon} = 0$ and $\div \vec v_\text{recon} = 0$ uniformly.
\end{theorem}
\begin{proof}
 Consider the biparabolic reconstruction of equation \eqref{reconSum2dRect} and Theorem \ref{reconRecThm} together with equation \eqref{eq:value9}. Applying the discrete Fourier transform one can rewrite
 \begin{align}
   q_{\text{recon}, ij}(\vec x)&=\sum_{m=1}^{9}c_{m,ij} b_m(\vec x)\\
   \begin{split} &= \left [ \hat q\no\left( \frac{1}{t_x t_y} \left(b_1(\vec x) +  \frac{3}{16} b_9(\vec x)  \right) 
   + \frac{1}{t_y} \left(b_3(\vec x) + \frac{3}{16} b_9(\vec x)  \right) \right . \right .\\
   \\ &\qquad\qquad + \left. b_5(\vec x) + \frac{3}{16} b_9(\vec x) 
   + \frac{1}{t_x} \left(b_7(\vec x)  + \frac{3}{16} b_9(\vec x)  \right) \right )\\
   &+ \hat q\eh \left( \frac{1}{t_y} \left(b_2(\vec x) -  \frac{3}{4} b_9(\vec x)  \right)  
   +  \left(b_6(\vec x)  - \frac{3}{4} b_9(\vec x)  \right) \right )\\
   &+ \hat q\ev \left(  \left(b_4(\vec x) -  \frac{3}{4} b_9(\vec x) \right ) 
   + \frac{1}{t_x} \left(b_8(\vec x) - \frac{3}{4} b_9(\vec x)  \right) \right )\\
	  & \left. +\frac{9}{4} b_9(\vec x)  \hat q\avg \right ] \cdot t_x^i t_y^j \end{split}
\end{align}
$p_\text{recon}(\vec x) = 0$ is thus clear. Moreover, one easily can compute $\div \vec v_\text{recon}(\vec x)$ by differentiating the basis functions $b_m$ as given in theorem \ref{reconRecThm} and verify $\div \vec v_\text{recon} = 0$. The computation is lengthy but uneventful, and is thus omitted. 
\end{proof}
\emph{Note:} From \eqref{eq:statinitialdatafourier4}, \eqref{eq:statinitialdatafourier1} it follows that
 \begin{align}
  \hat u\no \frac{(t_x-1)(t_y+1)}{\Delta x} + \hat v\no \frac{(t_x+1)(t_y-1)}{\Delta y} &= 0 \label{eq:discdivno}\\
  \hat u\avg \frac{1 + 4 t_y + t_y^2}{t_y} \frac{(t_x-1)(t_x+1)}{\Delta x t_x} + \hat v\avg \frac{1 + 4 t_x + t_x^2}{t_x} \frac{(t_y-1)(t_y+1)}{\Delta y t_y} &= 0
 \end{align}
 and from \eqref{eq:statinitialdatafourier2}--\eqref{eq:statinitialdatafourier3}
 \begin{align}
  \hat u\ev \frac{t_x-1}{\Delta x t_x} + \hat v\eh \frac{t_y-1}{\Delta y t_y} &= 0 \label{eq:discdiveh1}\\
  \hat u\eh \frac{1 + 6 t_y + t_y^2}{t_y} \frac{t_x-1}{\Delta x} + \hat v\ev \frac{1 + 6 t_x + t_x^2}{t_x} \frac{t_y-1}{\Delta y} &= 0 \label{eq:discdiveh2}
 \end{align}
 These equations are discretizations of $\div \vec v = 0$. 
 
\emph{Note 2:} Having thus found one eigenvector that belongs to the eigenvalue 0 of $\mathcal A(t_x, t_y) - \id_{12\times12}$ it is not a priori clear whether there are more. This can be checked by actually computing the kernel of $\mathcal A(t_x, t_y) - \id_{12\times12}$. Due to the extreme length of the expressions, the one-dimensionality of the kernel could so far only be verified using \textsc{mathematica}.

\begin{corollary}\label{thm:statpreservAF}
 The active flux scheme for linear acoustics in two spatial dimensions is stationarity preserving.
\end{corollary}
\begin{proof}
The proof consists of two parts. First, stationarity of the average values has to be checked. Second, the multi-step integration procedure (section \ref{chapterRectangularFlux}) has to be taken into account.
\begin{enumerate}[i)]
 \item Assume the initial data to fulfill \eqref{eq:statinitialdatafourier1}--\eqref{eq:statinitialdatafourier4}. By the above theorem this implies stationarity of the {point values at cell boundaries} and thus the three brackets in \eqref{eq:fluxbisimpson} are all equal. The values of $p$ are all zero, and thus the fluxes of $u$ and $v$ are zero. The change of the cell average of the pressure is
 \begin{align} \begin{split}
  &\left[ \frac16 \hat u\no \left( 1 + \frac{1}{t_y} \right ) + \frac{4}{6} \hat u\ev \right ] \frac{1}{\Delta x} \left( 1 - \frac{1}{t_x} \right ) \\  &\qquad\qquad +  \left[ \frac16 \hat v\no \left( 1 + \frac{1}{t_x} \right ) + \frac{4}{6} \hat u\eh \right ] \frac{1}{\Delta y} \left( 1 - \frac{1}{t_y} \right )  \end{split} \label{eq:fluxdifference}
 \end{align}
 Equations \eqref{eq:discdivno}, \eqref{eq:discdiveh1} imply that \eqref{eq:fluxdifference} vanishes identically. Thus the cell averages are stationary.
\item Assume again that the initial data fulfill \eqref{eq:statinitialdatafourier1}--\eqref{eq:statinitialdatafourier4}. Then at time step $n+\frac12$ and at time step $n+1$ they are equal to those at time step $n$. This is because they both are computed from the same initial data at time step $n$. The update of the cell average happens only at time step $n+1$. Therefore stationarity preservation is valid independently of how many steps are used for the integration in time.
\end{enumerate}
This completes the proof.
\end{proof}

Stationarity preservation has several consequences: The eigenvector $\hat Q$ of theorem \ref{thm:stationaryeigenvector} is the Fourier transform of those data that the active flux scheme keeps exactly stationary. In \cite{barsukow17a} it is shown that many numerical schemes add so much diffusion that only trivial (e.g. constant) stationary states are stationary at the numerical level as well. The active flux scheme on the other hand is stationarity preserving, i.e. it keeps stationary a discretization of all the stationary states of the PDE. By inverting the Fourier transform one obtains from \eqref{eq:discdivno}--\eqref{eq:discdiveh2} the following discrete relations that characterize the discrete stationary states:

\begin{align}
  \frac{ \{ [u\no]_{i+\frac12} \}_{j+\frac12} }{\Delta x} + \frac{ [ \{ v\no \}_{i+\frac12} ]_{j+\frac12} }{\Delta y} &= 0 &
  \frac{ \langle [ u\avg ]_{i\pm1} \rangle^{(4)}_{j}  }{\Delta x } + \frac{ [ \langle v\avg \rangle^{(4)}_{i} ]_{j\pm1} }{\Delta y} &= 0\\
  \frac{ \langle [   u\eh ]_{i+\frac12} \rangle^{(6)}_{j}  }{\Delta x} + \frac{ [  \langle v\ev \rangle^{(6)}_{i} ]_{j+\frac12} }{\Delta y} &= 0 &
  \frac{ [ u\ev ]_{i-\frac12,j} }{\Delta x } + \frac{  [v\eh_i ]_{j-\frac12}}{\Delta y } &= 0 
 \end{align}

 Here the following notation has been used:
 \begin{align}
  [q]_{i+\frac12} &= q_{i+1} - q_i &  \{q\}_{i+\frac12} &= q_{i+1} + q_i \\
  [q]_{i\pm1} &= q_{i+1} - q_{i-1} &  \langle q \rangle^{(\alpha)}_{i} &= q_{i-1} + \alpha q_i + q_{i+1}
 \end{align}
 The notation is combined for different directions, e.g.
 \begin{align}
  \{ [q]_{i+\frac12} \}_{j+\frac12} &= q_{i+1,j+1} - q_{i,j+1} + q_{i+1,j} - q_{ij}
 \end{align}

By \cite{barsukow17a} a stationarity preserving scheme is vorticity preserving. This means that there exists a discretization of $\nabla \times \vec v$ that is kept stationary by the scheme (even if the solution itself is not stationary). This is a discrete counterpart to Equation \eqref{eq:stationaryvorticity}. The discrete Fourier transform of the discrete vorticity is given by the left eigenvector $\Omega$ belonging to eigenvalue zero. As $\mathcal A$ depends on $\Delta t$, one is facing the slightly surprising situation that the discrete vorticity might depend on $\Delta t$. Indeed, this might have occurred for the discrete stationary states already, but turns out not to be the case. Additionally, the amount of computations made it impossible to determine the explicit shape of the discrete vorticity. However, its existence is clear by the results of \cite{barsukow17a}. From the dimension of the kernel it is also clear that there is only one discrete vorticity that is invariant. Only this particular discrete vorticity will be kept stationary while other discretizations will undergo some evolution. 

That the discrete vorticity depends on $\Delta t$ can be shown as follows: Consider an expansion of $\mathcal A$ in powers of $\Delta t$:
\begin{align}
 \mathcal A - \id = \Delta t \mathcal A^{(1)} + \Delta t^2 \mathcal A^{(2)}  + \Delta t^3 \mathcal A^{(3)}
\end{align}
Independently of $\Delta t$, the vector $\hat Q$ from theorem \ref{thm:stationaryeigenvector} is a right eigenvector of $\mathcal A - \id$ corresponding to an eigenvalue 0. Therefore one concludes that $\mathcal A^{(1)} \hat Q = \mathcal A^{(2)} \hat Q = \mathcal A^{(3)} \hat Q = 0$. It is possible to compute the left eigenvector $\Omega^{(1)}$ of $\mathcal A^{(1)}$:
\begin{align}
 \Omega^{(1)} &= \left( 0,0,0,   -\frac{t_y-1}{\Delta y t_y},0,0,     0,\frac{t_x-1}{\Delta x t_x},0,    0,0,0 \right )
\end{align}
This would correspond to a discrete vorticity
\begin{align}
 \frac{[v\ev]_{i-\frac12,j}}{\Delta x} -\frac{[u\eh_i]_{j-\frac12}}{\Delta y} \sim \del_x v - \del_y u
\end{align}

One can verify, however, that $\Omega^{(1)}  \mathcal A^{(2)}  \neq 0$ and thus $\Omega \neq \Omega^{(1)}$.

Finally, low Mach compliance is another consequence of stationarity preservation as discussed in \cite{barsukow17a}. Von Neumann stability of the method implies that no Fourier mode is growing. Thus the long time solution is governed by the discrete stationary states, and stationarity preserving schemes discretize all of the analytic stationary states (for more details see \cite{barsukow17a}). As $\epsilon$ goes to zero, one thus observes that the numerical solution obtained with the active flux scheme is a discretization of the limit solution. 

{In order for finite volume methods to be stable under explicit time integration upwinding is introduced. Many numerical schemes then strongly diffuse stationary states, and equivalently are unable to resolve the low Mach number limit. One remedy is to use implicit time integration without upwinding (central discretizations) (e.g. \cite{cordier12} and many others). If fully explicit time integration is required (e.g. if both low Mach number regions and shocks are present in the flow), then many numerical schemes are made low Mach compliant by introducing some kind of fix (e.g. \cite{dellacherie16,barsukow16}). For the active flux scheme presented here such a fix is not needed: the exact evolution operators at the boundary use the correct directions of information propagation while not spoiling the low Mach number limit.}

\section{Summary}

The active flux scheme is a finite volume method with additional degrees of freedom located on the cell boundary. Introduced in \cite{vanleer79} for one-dimensional linear advection, it has since been extended to triangular grids (\cite{eymann13}) and other systems of hyperbolic PDEs. This paper presents an implementation for two-dimensional Cartesian grids. This underlines the viewpoint that the active flux is a concept that can be used with considerable flexibility with respect to the computational grid.

Cartesian grids have several advantages. Apart from the ease of implementation and low computational cost, the limit of low Mach numbers can be studied very efficiently for linear acoustics on Cartesian grids (\cite{barsukow17a}). Low Mach compliance of a numerical scheme for linear acoustics is linked to it being \emph{stationarity preserving}, i.e. to be able to keep stationary a discretization of all the analytic stationary states. Here, the framework of \cite{barsukow17a} is used in order to show that the active flux possesses this property, and thus is low Mach compliant. The stationary states of the active flux scheme are obtained explicitly. They are closely linked to the solutions of the active flux scheme in the limit of low Mach numbers. 

Finally, as stationarity preservation is equivalent to vorticity preservation, it is thus possible to show that the active flux scheme for linear acoustics is vorticity preserving. This can often be checked easily if it is known which discretization of the vorticity remains stationary. The approach via stationarity preservation, on the other hand, allows to check vorticity preservation without having to know the discrete vorticity in advance. {We expect analogous results to hold true for the active flux scheme on three-dimensional Cartesian grids as well.}

Several numerical examples show the good performance of the active flux scheme in practice. In particular they show that active flux seems to be stable with the maximum CFL condition. To study this theoretically is subject of future work. As the scheme is of high order limiting is a further aspect of future investigation. 

Low Mach compliance is shown for the scheme endowed with an exact evolution operator for the {point values at cell boundaries}. This is different from Riemann solver based schemes, where even an exact Riemann solver does not in general lead to a low Mach compliant scheme. This, therefore, is clearly an advantage of the active flux scheme. However, the exact solution operator is generally unavailable for nonlinear systems of hyperbolic PDEs. Therefore approximate evolutions will be necessary. With the results of this paper in mind, future work shall focus on choosing them such that low Mach compliance is retained.

\newcommand{\etalchar}[1]{$^{#1}$}

\end{document}